\input amssym.def
\input amssym

\font\de=cmssi10

\def\pp{periodic point}
\def\po{periodic orbit}
\def\fo{finite-order}
\def\PNE{periodic Nielsen equivalent}
\def\SNE{strong Nielsen equivalent}
\def\NE{Nielsen equivalent}

\def\ra{\rightarrow}
\def\ie{{\it i.e.}}
\def\eg{{\it eg.}}
\def\cf{{\it cf.}}
\def\N{{\Cal N}}
\def\F{{\Cal F}}

\def\Na{N^\ast}
\def\ha{h^\ast}
\def\Ra{R^\ast}
\def\Nc{N^\circ}

\def\reals{{\Bbb R}}

\def\ne{\buildrel{\scriptscriptstyle NE}\over  \sim}
\def\pne{\buildrel{\scriptscriptstyle PN}\over  \sim}
\def\sne{\buildrel{\scriptscriptstyle SN}\over  \sim}

\def\circle{S^1}
\def\tf{\tilde{f}}
\def\th{\tilde{h}}
\def\tb{\tilde{b}}
\def\td{\tilde{d}}
\def\inv{^{-1}}
\def\tx{\tilde{x}}
\def\ty{\tilde{y}}
\def\tz{\tilde{z}}
\def\tM{\tilde{M}}
\def\tA{\tilde{A}}
\def\tN{\tilde{N}}
\def\tp{\tilde{p}}
\def\tq{\tilde{q}}
\def\tg{\tilde{g}}
\def\trho{\tilde{\rho}}
\def\tphi{\tilde{\phi}}
\def\tPhi{\tilde{\Phi}}
\def\tgamma{\tilde{\gamma}}
\def\tGamma{\tilde{\Gamma}}

\input epsf
\input amstex
\documentstyle{amsppt}
\NoBlackBoxes


\topmatter
\title Isotopy Stability of Dynamics on Surfaces\endtitle
\author Philip Boyland\endauthor
\leftheadtext{PHILIP BOYLAND}%

\address Department of Mathematics,
University of Florida,
Gainesville, FL 32611-8105\endaddress

\email boyland\@math.ufl.edu\endemail

\thanks  Thanks to Boju Jiang for pointing
out the necessity of part (c) in Lemma 1.3.\endthanks



\abstract
This paper investigates dynamics that persist
under isotopy in classes of orientation-preserving homeomorphisms
of orientable surfaces. The persistence of \pp s with
respect to periodic and strong Nielsen equivalence is studied.  
The existence of a dynamically minimal representative with respect
to these relations is proved and  necessary and sufficient
conditions for the isotopy stability of an equivalence class are given.
It is also shown that most the dynamics of the minimal representative
persist under isotopy in the sense that any isotopic map has
an invariant set that is semiconjugate to it.
\endabstract

\endtopmatter

\document

\head 0. Introduction\endhead
Isotopy stability of dynamics refers to dynamical behavior that persists
under isotopy. Since this behavior is present in 
every homeomorphism in an isotopy class, results about isotopy stability
often allow one to gain a great deal of dynamical information about a  map
given only fairly rough algebraic or combinatorial data about its isotopy
class.  For example, using the Lefschetz Fixed 
Point Theorem one can algebraically compute that every map in an
isotopy class has a fixed point.
Nielsen's  work on homeomorphisms of surfaces introduced what is
now called the Nielsen class of a fixed point. This work has been
generalized and is the content of Nielsen Fixed Point Theory.
The use of Nielsen classes
yields a refinement of the Lefschetz Theorem that often gives persistence of a
larger collection of fixed points.
Somewhat surprisingly, the application of this theory to \pp s is 
comparatively recent with the work of Halpern [Hp2], 
Jiang [J3], and others.  The proceedings [Mc] give a good sense of 
the current state of the theory.

For \pp s the analysis of isotopy stability  
begins with a definition of an equivalence relation.
The second step is to specify what it means for equivalence classes
of \pp s from isotopic maps to correspond. The basic persistence result
then gives conditions on an equivalence class which insure that any 
isotopic map has a nonempty equivalence class of \pp s that 
corresponds to the given one. This makes precise the meaning of
``present in every element of the isotopy class'', and 
it allows one to consider certain classes as invariants of isotopy 
or as part of the isotopically stable dynamics.
Given such a persistence result, the next question is that of a dynamically
minimal representatives; is there a map in each isotopy class
that has only the forced dynamics and nothing more? More specifically, is
there a map which has exactly one \pp\ in each isotopically stable
equivalence class and no other periodic points?

The first two sections of this paper concerns theories of 
this type generated by the two equivalence relations 
of periodic and strong Nielsen equivalence.
Results on the first of these theories are contained in
[J3], [HPY], and  [HY] and on the second in  [AF] and [Hll]. 
In the case of homeomorphisms of surfaces
this paper strengthens and extends certain of these results. 
In particular, the existence
of a dynamically minimal representative is proved as well as 
necessary and sufficient conditions for an equivalence class of \pp s
to be isotopically stable.

The dynamically minimal representative in an isotopy class is
a refinement of the Thurston-Nielsen canonical form.
There have been a number of papers that have refined this canonical
form for dynamical purposes, for example, [S],  [BS], and [BL].
A refinement of the Thurston-Nielsen canonical form was also
used to prove the existence of a dynamically minimal 
representative for Nielsen classes of fixed points in the category
of surface homeomorphisms. This result
was sketched in [J2] and [I], and given in full detail in [JG].

The existence of a dynamically minimal representative
with respect to \pp s is somewhat more delicate than that of
fixed points. It requires that a single map have the minimal number 
of \pp s of all periods. The difficulty is that the map that has the 
least number of fixed points in the isotopy class of $f^n$ may
not be itself the n$^{th}$ iterate of a map that has the least number
of fixed points in the class of $f$.

As a simple example, let $\phi$ be a map
of a surface with $\phi^n = Id$ for some (least) $n$.  Further suppose 
that $\phi$ has several \po s that have 
period less than $n$. As a consequence of Lemma 1.1(aii) and Lemma 2.3 these
\pp s are isotopically stable  and thus must be present
in any dynamically minimal model.  On the other hand, since
$\phi^n$ is the identity there is a map in the isotopy class of $\phi^n$
with just one fixed point. In this example there is {\it no}  map $g$ isotopic
to $\phi$ with the
property that $g^n$ has the least number of fixed points in its class  for
all $n$. However, Theorem 2.4 shows that there is a homeomorphism isotopic 
to $\phi$ that has the least number of \pp s of each period.

The last section discusses a kind of uniformization of the persistent 
\pp s. These results extend those
of [H] and [Ft]  to reducible mapping classes. Roughly speaking, one
obtains the persistence of orbits that are not periodic
by taking the closure of the isotopically stable \pp s. 
This yields the isotopy stability of essentially all the nontrivial
dynamics of the dynamically minimal representative.

The first results on isotopy stability 
that have this kind of global character appear to be Franks' work on Anosov
diffeomorphisms and Shub's work on expanding maps. R. MacKay pointed
out to the author that within Differential
Geometry results of this type are much older.
They go back at least to a paper of Morse on geodesics
on surfaces ([M]). A survey of dynamical applications of isotopy stability
is given in [Bd2].

\head 1.  Condensed homeomorphisms\endhead
This section develops a refinement of the Thurston-Nielsen
canonical form for isotopy classes of maps on surfaces. There
are two steps in the refinement. The first step involves
an alteration of the behavior of the map on the closed reducing
annuli between components. The new map is 
called an adjusted reducible map and is described in
Lemma 1.3. The second step is more radical and replaces
the map on a \fo\ component by a dynamically simpler map and may involve
slightly altering the topology of pseudo-Anosov components
and identifying points in adjacent components. The result
of this process is called a  condensed homeomorphism and
is described in Theorem 1.7.

A basic familiarity with Thurston's work on surfaces 
is assumed. For more  information see
[T], [FLP], or [CB]. For more information on Nielsen fixed
point theory and the analogous theory for \pp s, see [J3].

Throughout this paper $M$ will be a compact, orientable
 $2$-manifold with negative Euler characteristic and  perhaps
with boundary. Unless otherwise
noted, all homeomorphisms $M\ra M$ will
be orientation-preserving, and 
isotopies do not need to fix the boundary
of $M$ pointwise. The {\de period} of a \pp\ always means its
{\it least \/} period, and 
the notation $per(x, f) = n$ indicates that $x$ is a \pp\ with 
period $n$ of a  homeomorphism $f:M\ra M$.  The set of all \pp s with period
$n$ is denoted $P_n(f)$, and  the set of points fixed by $f$ is
$Fix(f)$. Note that, in general, $Fix(f^n)$ may be larger than
$P_n(f)$.

 Given  $x, y \in Fix(f)$, $x$ is {\de Nielsen equivalent} 
to $y$ (written $(x,f)\ne (y,f)$,
or if the map is clear from the context, $x\ne y$) if there
is an arc $\gamma:[0,1] \ra M$ with $\gamma(0) = x$, $\gamma(1)
=y$, and $f(\gamma)$ is homotopic to $\gamma$ with fixed
endpoints . In this case we say that $x$ is {\de  \NE\ to\/} $y$ {\de via\/}
$\gamma$. 
Given  $x,y\in P_n(f)$, $x$ is {\de \PNE\ } to $y$
if $(x,f^n) \ne (y,f^n)$ (written $(x,f)\pne (y,f)$).
It is important to note that our definition of
 periodic Nielsen equivalent requires that the \pp s have 
the same least period. This requirement is not completely
standard in the literature. Two periodic {\it orbits\/} 
are \PNE\ if \pp s from each orbit are.
In what follows, equivalence of both \pp s and \po s will be considered.
A certain amount of confusion will probably be avoided if the distinction
between these two notions is maintained. Also note that here we
primarily consider the geometric notion of Nielsen classes for \pp s.
In particular, the phrase ``periodic Nielsen class'' always refers
to a {\it nonempty} equivalence class of \pp s.
For a general account of the algebraic theory see [J3], [HPY], and  [HY].

If $x$ and $y$ are \pp s with $(x,f^n) \ne (y,f^n)$
but  $n = per(x,f) > per(y,f)$, $x$ is said to {\de collapse} to
$y$ (written $(x,f) \vdash (y, f)$ or  $x \vdash y$). 
 One periodic {\it orbit\/} is collapsible to another
if \pp s from each orbit are.

If the Nielsen equivalence in the various definitions is
realized by an arc $\gamma$, this is indicated by
saying that $x \pne y$ or $x \vdash y$ {\de via the arc} $\gamma$.
The \pp s $x$ and $y$ are said to be  {\de related} if 
$x\pne y$, $x\vdash y$, or $y\vdash x$ .
A simple consequence of the definition is that when   $x$ and
$y$ are related,    $(x, f^k)\ne (y, f^k)$ for  
$k=\max\{per(x,f), per(y,f)\}$. 

It will be useful to extend the notion of periodic Nielsen
equivalence to certain closed curves in $M$. Assume that the
closed curve
$C$ is fixed setwise by $f$ and $x$ is a fixed point of $f$. The
{\de invariant curve}  $C$ {\de is Nielsen equivalent to} $x$ if there is an
arc $\gamma:[0,1]\ra M$ with $\gamma(0) \in C$, $\gamma(1) =
x$, and $f(\gamma)$ is homotopic to $\gamma$ via a homotopy $F_t$
with $F_t(0)\in C$ and $F_t(1)=x$ for all $t$. The notation for this
is  $(x,f)\ne (C, f)$. The notions of periodic
Nielsen equivalence between a point and a curve, equivalence of
two periodic curves, and  the analogous notions of collapsible are defined
in the obvious way.

If $\phi:M\ra M$ is an isometry of a hyperbolic metric, then it is
standard that $\phi$
is {\de finite-order}, {\it i.e.} there is some least $n>0$ (called
the {\de period}) with
$\phi^n = Id$. Conversely, when $\phi$ is \fo, it is topologically
conjugate to an isometry
of some hyperbolic metric. In the literature  \fo\
homeomorphisms  are often called ``periodic'', but that terminology
is avoided for obvious reasons. A \pp\ of a \fo\
homeomorphism that has the
same period as the homeomorphism is called {\de regular}; any
point with lesser period is called a {\de branch} \pp.
Since only orientation-preserving \fo\ maps are considered
here, the set of branch \pp s is always finite.

A homeomorphism  $\phi$ is called {\de pseudoAnosov (pA)} if there
exists a number $\lambda > 1$ (called the {\de expansion constant})
and a pair of transverse measured
foliations $(\F^u, \mu^u)$ and $(\F^s, \mu^s)$ with
$\phi(\F^s, \mu^s) = (\F^s, {1 \over \lambda}\mu^s)$ and
$\phi(\F^u, \mu^u) = (\F^u,  \lambda\mu^u)$. 
For a pA map on a boundary component there is a certain
amount of choice involved in the structure of the foliation 
and with the dynamics.
Before giving the choices made here we describe
what will be called the {\de standard model\/} (the author learned
this succinct description from D. Fried). 

Roughly speaking, 
in a neighborhood of a boundary component the standard model of a pA
map looks like the map and foliation 
obtained by blowing up a singularity of a measured foliation.
This description is only informal because the blown up map will be continuous
only when the original map is differentiable 
at the singularity.  This is not the case in the usual
constructions  of pA maps (but {\it cf.\/} [GK]).
The main observation needed to avoid this difficulty is that
a pA map is differentiable at regular points and 
the behavior at other singularities can be obtained by using
branched covers.

There are three main cases for boundary behavior that
correspond to three situations for an
interior singularity: the singularity is fixed by the pA map
and the prongs do not rotate, the singularity is fixed and the
prongs rotate, and the singularity is a \pp.
We first describe the standard model on the boundary that corresponds
to the case of a fixed, nonrotating prong.

Let $L:\reals^2\ra\reals^2$ be the linear map with matrix
$$\pmatrix \lambda&0\\ 0&{1\over\lambda}\endpmatrix $$
Replace the origin of $\reals^2$ with a boundary circle 
$C$  and let $R^\prime$ be the resulting space.
Now use the action of $L$ on lines to 
define a homeomorphism $L^\prime:R^\prime\ra R^\prime$.
The stable and unstable foliation of $L$ induce a stable and unstable
foliation on $R^\prime$.
These stable and unstable foliations and the homeomorphism $L^\prime$
give the standard model for a pA map in the neighborhood of 
a boundary component that corresponds to a regular point.
For a boundary component that corresponds to a one-prong,
take this model and project it via $z\mapsto z^2$. 
For a boundary component that corresponds to a $k$-prong,
take this one-prong model and lift it to 
the cover that has projection  $z\mapsto z^k$.
Figure 1a shows the dynamics at an interior three-prong singularity and
Figure 1b shows  an analogous boundary
component in the standard model. Note that the curves in
these figures indicate the motion of orbits and not necessarily leaves
of the foliations. The segments between the
singularities on the boundary will be termed  {\de degenerate
leaves} and are considered leaves of both the stable and
unstable foliation.  One has to allow for this degeneracy in
the definition of {\it transverse\/} measured foliations. 

\topinsert
\def\epsfsize#1#2{.7\hsize}
\centerline{\epsfbox{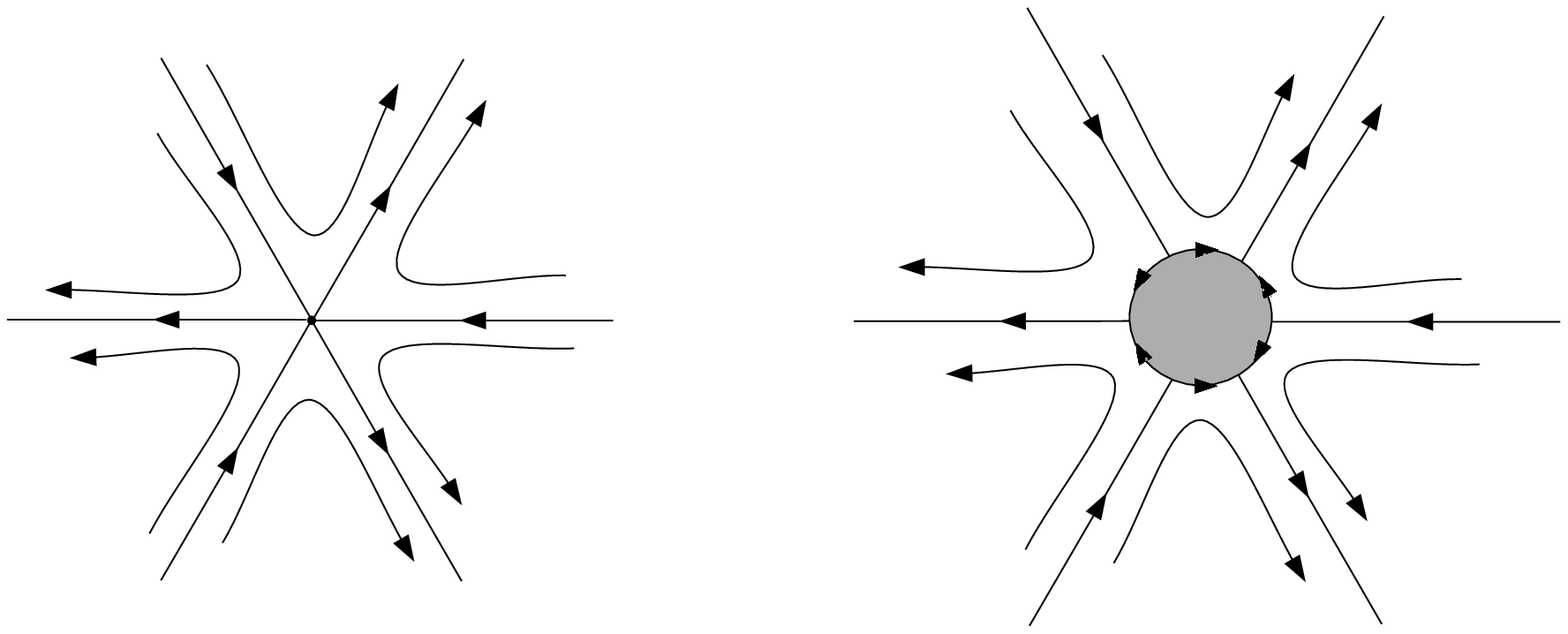}}
\botcaption{Figure 1}{{\bf (a)} Dynamics of a pA map at a three-prong singularity.
{\bf (b)} Dynamics at the boundary analog of (a) in the standard model
of a pA map.}
\endcaption
\endinsert

\topinsert
\def\epsfsize#1#2{.7\hsize}
\centerline{\epsfbox{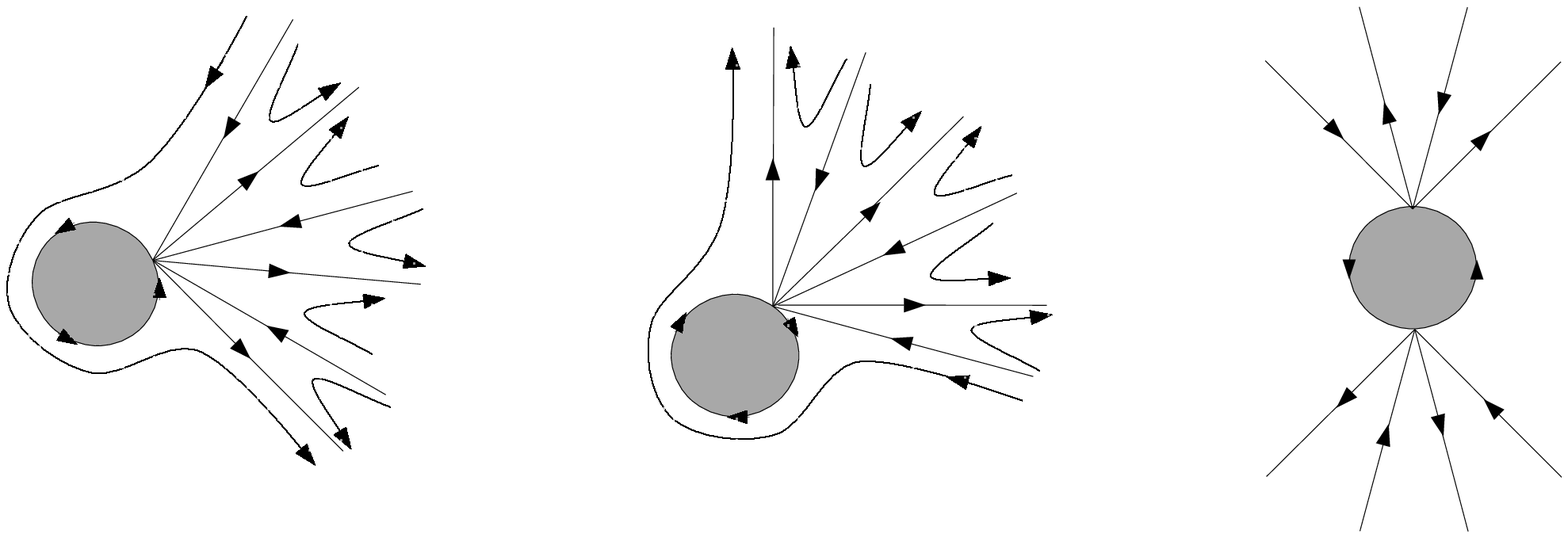}}
\botcaption{Figure 2}{{\bf (a)} The dynamics of a boundary-adjusted pA map at
the boundary analog of a non-rotating three-pronged singularity.
{\bf (b)} Same as (a) but the result of a different collapse.
{\bf (c)} The dynamics of a boundary-adjusted pA map at
the boundary analog of a four-pronged singularity with rotation number
$1/2$.}
\endcaption
\endinsert

As constructed, in a boundary component of the standard model that
corresponds to a fixed, nonrotating
singularity,  the degenerate leaves consist of orbits that are
heteroclinic from the intersection of the stable leaves with
the boundary to the intersection of the unstable leaves with the
boundary.
(One may also adjust things to make the map the identity on
these segments, see [GK] and [JG].) In the analog of fixed,
rotating singularity, this model on the boundary is composed
with the appropriate rigid rotation. In the analog of  a
singularity that is a \pp, this model is composed with the
appropriate translation of one boundary component onto another.

Since a goal of this paper is to create models in each isotopy
class that have the least dynamics possible, the standard model
on the boundary needs some minor adjustments.
For the analog of
a fixed, nonrotating prong the
necessary modification is achieved by  collapsing down all the
degenerate leaves except one. The dynamics in this case are
just those that descend from the standard model, \ie,  the boundary
consists of a fixed point and  a homoclinic loop. There are two
non-conjugate choices for the collapse. These are shown in Figures 2(a) and
2(b) for the boundary analog of a 
non-rotating three-prong singularity. In the analog of a fixed,
rotating $m-$prong the form of the collapse depends on the
rotation number of the orbits on the boundary circle in the standard
model. If this rotation number is $p/q$ with $p$ and $q$ relatively
prime,  one collapses down collections of adjacent groups of $2 m/q -1$
degenerate leaves. This will leave a single period $q$ orbit on the
boundary whose points are connected by homoclinic segments.
The result of the
collapse for a
four-prong singularity  with rotation number $1/2$ is shown in Figure 2(c).

If the
boundary component is moved off itself by the pA map, the appropriate
collapse is chosen based on whether the prongs are rotated or
not when the component's forward orbit first lands on itself.
PseudoAnosov homeomorphisms with this collapsed boundary behavior
will be called {\de boundary-adjusted} pA. It is clear that a
boundary-adjusted pA and the corresponding standard model are
conjugate on the interior of $M$.

The first lemma describes the relations among the \pp s of
pA and \fo\ maps.

\proclaim{Lemma 1.1} 

(a) If $\phi: M\ra M$ is \fo, then

{\leftskip=60pt\rightskip=20pt\parindent=-18pt

(i) Each regular \pp\ is \PNE\ to every other regular periodic
point.

(ii) All branch \pp s are unrelated to each other and are not
\PNE\ or collapsible to any boundary component.

(iii) Each regular \pp\ is collapsible to any branch \pp.

}
\smallskip
(b)  If $\phi: M\ra M$ is boundary-adjusted pA, then

{\leftskip=60pt\rightskip=20pt\parindent=-18pt

(i) Each interior  \pp  is  unrelated to any other \pp.

(ii) Each boundary \pp\ is \PNE\ to every \pp\ on its orbit that
is contained in the same boundary component and is unrelated to any
other \pp.

(iii) Each boundary component is unrelated to any \pp\ except the \pp s
it contains and is unrelated to any other boundary component.

(iv) If a boundary component $b$ is \PNE\ to itself via an arc
$\gamma$, then $\gamma$ is null homotopic via a homotopy that
constrains the endpoints of $\gamma$ to lie on $b$.

}
\endproclaim

\demo{Proof} (a) ({\it cf.} [J1], Section  7) We prove (ii) first. Assume
$x$ and $y$ are two branch periodic orbits and $(x, \phi^k) \ne
(y, \phi^k)$ via an arc $\gamma$ with $k$ less than
the period of $\phi$. Since $\phi$ is a hyperbolic
isometry, the unique geodesic isotopic to $\gamma$ with fixed
endpoints is fixed by $\phi^k$. Because $\phi$ is
orientation-preserving, this implies that $\phi^k = Id$,
a contradiction.
The second part of (ii) proved like the first after  blowing down
the boundary components of $M$ to points. Parts (i) and (iii) are
easy consequences of the fact that $\phi^n = Id$.

(b) This has been essentially proven by many authors, {\it eg.}
[BK], [H],  [I], and  [JG]. We will remark in Section 3 how the
methods of [H] can be adapted to deal with the assertions
involving the boundary. 
\quad\qed
\enddemo

The Thurston-Nielsen classification theorem for isotopy classes
of surface homeomorphisms gives a (fairly) canonical
representative in each isotopy class. These representatives are
constructed from pA and \fo\ pieces glued together along annuli in which
twisting may occur. More precisely, a homeomorphism  $\phi$ is called {\de
reducible} if there exists  a collection of pairwise
disjoint simple closed curves, $\Gamma = \{\Gamma_1, \Gamma_2,
\dots, \Gamma_k\}$, in $int(M)$ with $\phi(\Gamma) = \Gamma$ and
each connected component of $M - \Gamma$ has negative Euler
characteristic. Further, $\Gamma$  comes equipped with a
$\phi$-invariant open tubular neighborhood $\N(\Gamma)$. The
connected components of $M-\N(\Gamma)$ are called the {\de
components} of $\phi$. The orbit of a component under $\phi$ is
called a $\phi${\de -component}. For a subset $X\subset M$, its
orbit under $\phi$ is denoted $o(X)$.

\proclaim{Thurston-Nielsen Classification Theorem} 
 Every orientation preserving homeomorphism of an orientable surface 
with negative Euler characteristic is
isotopic to a homeomorphism $\phi$ such that either

{\leftskip=40pt\parindent=-18pt

(a) $\phi$ is pseudoAnosov, or

(b) $\phi$ is  finite-order, or

(c) $\phi$ is reducible and on each $\phi$-component, $\phi$
satisfies (a) or (b).

}
\endproclaim

A map $\phi$  as in part (c) of the theorem above is called  
{\de TN-reducible}.  
The curves of $\Gamma$ are called {\de reducing curves}. The
connected components of $Cl(\N(\Gamma))$ are called {\de
reducing
annuli} and the reducing annulus containing $\Gamma_i$ is
denoted $A(\Gamma_i)$. 
Two components are called {\de adjacent}
if they each share a boundary curve with the same reducing
annuli. Note that a component can be adjacent to itself.

If $\Gamma_i$ is a reducing curve and for some $n$,
$\phi^n(\Gamma_i) = \Gamma_i$ with the  orientation on $\Gamma_i$ reversed,
then $\Gamma_i$ is called a {\de flipped reducing curve}. In this case
the reducing annulus $A(\Gamma_i)$ is called a {\de flipped reducing
annulus}. 

The next lemma provides tools for working with flipped reducing annuli.
The annulus $A= \circle\times [-1, 1]$ has universal cover
$\tA = \reals\times [-1, 1]$. For a orientation preserving
circle homeomorphism $f$, $\rho(f)$ denotes the rotation number of $f$.

\proclaim{Lemma 1.2} Let $f:A\ra A$ be an 
orientation preserving homeomorphism
which interchanges the boundary components of $A$.

{\leftskip=40pt\parindent=-18pt

(a) If for $i = 1, 2$, $b_i$ denotes $f^2$ restricted to
$\circle\times \{i\}$, then $\rho(b_{-1}) = - \rho(b_1)$

(b) There exists a  homeomorphism $g$ that is isotopic to $f$ rel 
the boundary of $A$, is equal to $f$ on the boundary of $A$, and
the only periodic points of $g$ in the
interior of $A$ are two fixed points that have nonzero index and
are not Nielsen equivalent.

}
\endproclaim

\demo{Proof}  Pick a lift $\tf:\tA\ra \tA$ and let 
$S:\tA\ra \tA$ be given by $S(x,y) = (-x, -y)$. Note that
$S^2 = Id$ and if $\th := S \tf$, then $\th$ is the lift of
an orientation and  boundary preserving homeomorphism
of $A$ and $\tf = S \th$. 
For $i = 1, 2$, let   $\th_i$ and $\tb_i$ denote
$\th$ and $\tf^2$ restricted to $\circle\times \{i\}$, respectively,
and $T::\reals\ra\reals$ be given by  $T(x) = - x$.
We then have $\tb_{-1} = T \th_1 T \th_{-1}$ 
and $\tb_1 = T \th_{-1} T \th_1$ which implies $(T h_1)^{-1} \tb_{-1} (T h_1)
= \tb_1$. Since $T \th_1$ is orientation reversing, 
$\rho(b_{-1}) = - \rho(b_1)$, proving (a).

To prove (b), let $S$ also denote be the projection of the involution
to $A$ and $h:= S\circ f$. In addition,
let $C_r = S^1 \times \{r\}$ for $r\in [-1,1]$, and for $i=-1,1$
let $\alpha_i$ be $h$ restricted to $C_i$. Since each $\alpha_i$ is
an orientation preserving circle homeomorphism, there is a
family of such homeomorphisms $\alpha_r$ interpolating $\alpha_{-1}$
to $\alpha_1$ so that
$\alpha_r = id$ for $r \in [-1/4, 1/4]$. Now define  $\alpha: A \rightarrow
A$  so that $\alpha$ restricted to each  $C_r$  is $\alpha_r$.  Thus
$h\circ \alpha^{-1} = id$ on the boundary of $A$, and so
by  the Alexander trick, $h\circ \alpha^{-1}\simeq id$ rel
the boundary of $A$.

Now let $\psi_t$ be flow shown in Figure 3 which we assume
has been constructed so that it commutes  with the involution $S$ and
is the identity on the boundary of $A$. If $g' = \psi_1\circ\alpha$,
 then $g' = \psi_1$ on $S^1 \times [-1/4, 1/4]$, and  any points not
in that sub-annulus and not on the boundary of $A$
will have their $\omega$-limit set under $g'$ in the sub-annulus.
Further, rel the boundary,
$g := S\circ g' \simeq S\circ\alpha \simeq S\circ h = f$ and
so $g$ is as required for part (b).
\quad\qed
\enddemo

\topinsert
\def\epsfsize#1#2{.7\hsize}
\centerline{\epsfbox{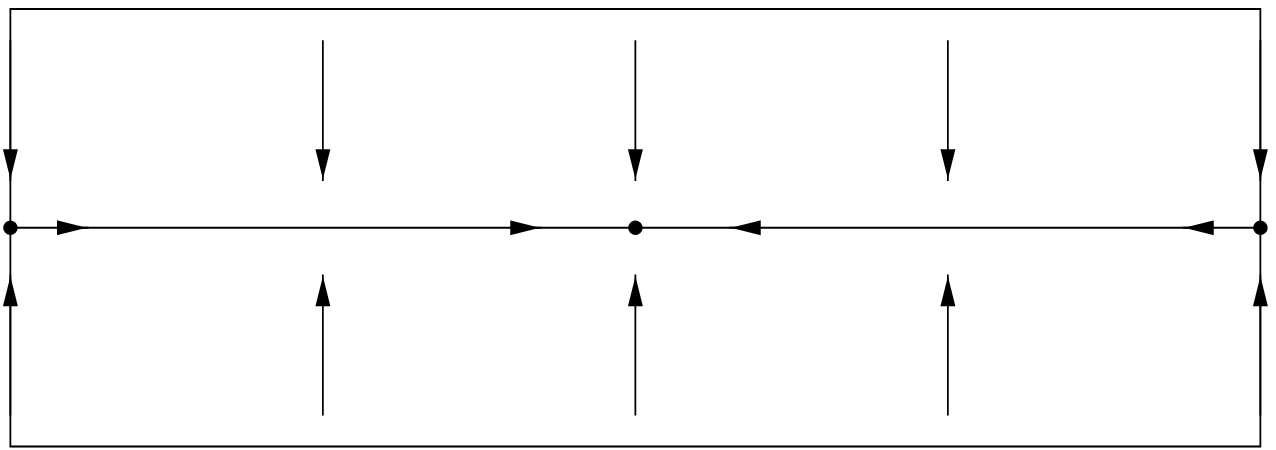}}
\botcaption{Figure 3}{The flow $\psi_t$ on the annulus
used in the proof of Lemma 1.2.
The fixed point at the left is at the point $(0,0)$.}
\endcaption
\endinsert

If $\phi$ is a TN-reducible map, a reducing annulus always has
\pp s on both of its boundaries. These \pp s are termed
{\de peripheral}. Note that peripheral \pp s for a \fo\ component
are always regular. If there are \pp s on distinct boundary
components of
a reducing annulus that are \PNE, the reducing annulus is called
{\de untwisted}.  An unflipped reducing annulus $A$ with
$\phi^n(A) = A$ is untwisted if and only if the boundary
components have the same rotation number in the lift of the
restriction of $\phi^n$ to the universal cover of the annulus.
A flipped reducing annulus $A$ is untwisted
if and only if the same condition holds for $\phi^{2n}$. Using Lemma 1.2
this  happens if and only if  $\phi^{2n}$ has fixed points on both
boundaries of $A$.

The next lemma gives the form of our first refinement of a
TN-reducible map.

\proclaim{Lemma 1.3} Any TN-reducible map, $\phi$, is isotopic to a reducible
map that satisfies:

{\leftskip=40pt\parindent=-18pt

(a) All pA components are boundary-adjusted.

(b) There are no \pp s in the interior of unflipped reducing annuli.

(c)  The interior of each flipped reducing annulus 
contains exactly two periodic points. These \pp s have nonzero index 
and are not periodic Nielsen equivalent.

(d) There are no untwisted reducing annuli connecting adjacent
\fo\ components.

}
\endproclaim

\demo{Proof} Starting with the TN-reducible map, $\phi$,
first replace all the maps on pA components by 
the appropriate boundary-adjusted pA. The construction of boundary-adjusted pA
maps makes it clear that this can be done within the isotopy class.
The behavior in part (b) can then be arranged in each unflipped
reducing annuli
by composing with a sufficiently strong push from one boundary
circle towards the other.
The behavior in part (c) can be arranged in each flipped
reducing annulus using Lemma 1.2(b). The last step is to show that
part (d) holds after eliminating unnecessary reducing annuli.

Let $\phi_1$ and $\phi_2$ be \fo\ maps 
on adjacent $\phi$-components $N_1$ and $N_2$ 
which each share a boundary with the
untwisted  reducing annulus $A$.  Because peripheral \pp s of
\fo\ components are
always regular for the component, the definition of untwisted
reducing annulus yields that $\phi_1$ and $\phi_2$ have the same
period, say $n$. Again using the fact that $A$ is untwisted,
$\phi$ restricted to  $N_1 \cup N_2 \cup o(A)$ must have its $n^{th}$ iterate
isotopic to the identity. A theorem stated by Nielsen (now contained in
the  classification theorem) implies that $\phi$ restricted to  
$N_1 \cup N_2 \cup o(A)$ is isotopic to a single \fo\ map of period $n$. Thus
the reducing annulus $A$ can be eliminated.  
\quad\qed
\enddemo

A reducible map as in Lemma 1.3  will be called  {\de adjusted}. 
If $\gamma$ is an arc connecting points $x$ and $y$,
$\gamma$ is said to {\de essentially intersect} a set $X$ if any
arc homotopic to $\gamma$ with fixed endpoints has nonempty
intersection with $X$.

\proclaim{Lemma 1.4} Let $\phi_1$ be a adjusted
reducible
map. If $x$ and $y$ are periodic points of $\phi_1$  that are
not contained in the interior of a flipped reducing annulus
and either $x\pne y$ or $x \vdash y$ via the arc $\gamma$ then:

{\leftskip=40pt\parindent=-18pt

(a) If $\gamma$ essentially intersects a reducing curve
$\Gamma_0$, then both $x$ and $y$ are related to $\Gamma_0$.
If $\gamma$ essentially intersects  two reducing curves
$\Gamma_1$  and $\Gamma_2$, then $\Gamma_1$  is related to  $\Gamma_2$.

(b) The arc $\gamma$ does not essentially intersect the interior
of any pA component.

(c)  If $\gamma$ essentially intersects a reducing curve
$\Gamma_0$, then the reducing annulus $A(\Gamma_0)$ is untwisted.

(d) The arc $\gamma$ cannot essentially intersect two \fo\
components. 

}
\endproclaim

\demo{Proof}
(a) This will be proved (independently) in Section 3 as the proof 
uses the covering space techniques introduced there.

(b)  Denote by $N$ the pA component that is essentially
intersected by $\gamma$.
There are two main  cases to consider. First assume that $x\in
Int(N)$. If $y$ is also in $Int(N)$, then Lemma 1.1b(i)  implies
that $\gamma$ must essentially intersect some reducing curve.
This must also be the case if $y\not\in Int(N)$. Let  $\Gamma_1$
denote the first such curve encountered as one traverses
$\gamma$ outward from $x$. If $b$ is a boundary component of
$N$ that is shared with $A(\Gamma_1)$, then (a) shows that
$x$ is related to $b$ contradicting Lemma 1.1b(iii).

The second case to be considered is when neither $x$ nor $y$ are
in $Int(N)$. In this case we can isotope $\gamma$ with fixed
endpoints to $\gamma_1$ that  essentially intersects two
boundary components of $N$, $b_1$ and $b_2$, with the portion of
$\gamma_1$ between the intersections contained in the interior
of $N$. If  $b_1 \not = b_2$,
then again using (a)  we get a
contradiction to Lemma 1.1b(iii) and if $b_1  = b_2$, a contradiction
to Lemma 1.1b(iv).

(c) By (b), the \pp s $x$ and $y$ can only be peripheral pA
or \fo. Thus if $\gamma$ does not essentially intersect a \fo\ component,
then $x$ and $y$ must be on the two boundary components of the
reducing annulus $A(\Gamma_0)$. Since for an annulus homeomorphism
\pp s on different boundary components cannot collapse to each other, the
only possibility is that $x\pne y$, and so $A(\Gamma_0)$ is untwisted.

Assume now that $\gamma$ is a simple arc and
essentially intersects a nonempty collection
of \fo\ components, $N_1, N_2, \dots, N_k$. Further assume that $y$ is
contained in a \fo\ component and 
$x$ is a peripheral pA \pp\ that lies on a closed curve $C_1$ that is the
boundary of a reducing annulus $A_1$. Let  $m$ be such that
$\phi_1^m(x) = x$ and $\phi_1^m = Id$ on all $N_i$. Now pick
$y^\prime$ contained on a boundary component $C_2$ of the \fo\ component which
contains $y$. There is an arc $\gamma^\prime$ so that
$(x, \phi_1^m) \ne (y', \phi_1^m)$ via $\gamma^\prime$, and $\gamma^\prime$
essentially intersects the same reducing curves as $\gamma$.

Now choose a pair of pants decomposition of $A_1\cup (\cup_{i=1}^k N_i)$
which refines the collection of reducing curves (with the exception of
the curve adjacent to $C_1$) and consider the Dehn-Thurston
parameterization of closed curves and  simple arcs with endpoints on
the boundary  using  this
decomposition (see [FLP], expos\'e 4 or [PH], \S 1.2) If any reducing curve that
essentially intersects $\gamma$ were contained in an twisted reducing
annulus, then the twist parameter of the parameterization at this
curve for $\gamma'$ and $\phi^m(\gamma')$ would be different. This
would imply that $\gamma'$ and $\phi^m(\gamma')$ are not homotopic, a
contradiction.

The proof when $x$ and $y$ are in other positions is similar. 
If the given $\gamma$ is not simple, pass to a finite cover
in which it is homotopic with fixed endpoints to a simple arc.
The result on the simple arc in this cover implies the result in the base.

(d) This  follows from (b), (c) and and the property
given in Lemma 1.3(d). 
\quad\qed
\enddemo

The next proposition describes all the relations among \pp s
of adjusted reducible homeomorphisms. A given set of
relations
among \pp s is said to {\de generate}  a (perhaps) larger
collection 
if the given set and the  following properties give
rise to all the relations.
\medskip

{\leftskip=40pt\parindent=-18pt

(1) Periodic Nielsen equivalence is an equivalence relation on
periodic points of the same period.

(2) $x\pne y$ and $y\vdash z$ implies $x \vdash z$.

(3) $x\vdash y$ and $y\pne z$ implies $x\vdash z$.

}

\medskip
Although it will not be needed here, it is perhaps worth noting
that these properties show that $\vdash$ induces an order
relation on
the set of periodic Nielsen equivalence classes. It is easy to see
that it is transitive and a minor alteration in the definition of
$\vdash$ (allow $per(x) \ge per(y)$) will give a partial order
on these equivalence classes. This is essentially the 
direct system of weighted sets discussed in [J3], Section III.3.

\proclaim{Proposition 1.5:} When the following relations exist
they generate all relations among \pp s of an adjusted
reducible homeomorphism $\phi_1$.

{\leftskip=40pt\parindent=-18pt

(a) A peripheral pA  \pp\  is \PNE\ to other \pp s on its orbit or 
to an adjacent peripheral pA or to an adjacent peripheral \fo\  \pp.

(b) A peripheral pA  \pp\  is collapsible to a periodic point
with half its period in the interior of an adjacent flipped, untwisted
reducing annulus.

(c) A regular \fo\ \pp\ is \PNE\  to another regular \fo\ \pp\ in the
same component.

(d) A regular \fo\ \pp\ is collapsible  to a branch \fo\
\pp\ in the same component.

}
\endproclaim

\demo{Proof} By Lemma 1.4(b) a \pp\ in the interior of a pA
component is unrelated to any other \pp. The remaining types of
\pp s are \fo,   peripheral pA and \pp s in the interior of
flipped reducing annuli. By Lemma 1.3(c), a \pp\ in the interior of
a flipped reducing annuli $A$ is not \PNE\ to the other \pp\ in the
interior of the same reducing annulus and is only related
to \pp s on the boundary  if $A$ is untwisted. Since $A$ is flipped
the adjacent components must have the same type. By Lemma 1.3(d),
if $A$ is untwisted this type cannot be \fo\ and the adjacent
components most therefore both be pA. As a consequence, (b) gives
the only possible relations for \pp s in the interior of
flipped reducing annuli.

Now by Lemma 1.4 (b) and (d), \fo\ and peripheral pA 
\pp s can only be related by a $\gamma$ that essentially
intersects the interior of at most one component and this component  must be
\fo.  Using Lemma 1.4(b) and (c), if a peripheral pA \pp\ $x$ is 
related to a \pp\ not on its orbit, $x$ must lie on the boundary of a
untwisted  reducing annulus and thus be periodic Nielsen
equivalent to a periodic point, say $z$,  on the other boundary
of $A$.  The case when the reducing annulus is flipped was
dealt with in the previous paragraph, so assume that the reducing annulus is
unflipped. If  the adjacent component is pA, again by Lemma 1.4(b),
$x$ is periodic Nielsen equivalent to $z$ and 
is unrelated to any \pp s not on the orbit of $x$ and $z$.
The second case is when the  adjacent component is \fo.
Now from  Lemma 1.4(b) and (d), a finite-order \pp\ is related
only to other \fo\ \pp s in the same component or to peripheral
\pp s in adjacent pA components.  Thus the only new relations
possible for $z$ (and thus for $x$ ) are: collapse to a branch
\pp\ in the same \fo\ component, periodic Nielsen equivalence to
a regular \pp\ in the same \fo\ component, and periodic Nielsen
equivalence to a peripheral \pp\ in an adjacent pA component.
All these relations are generated by those given in
the statement of the lemma. 
\quad\qed
\enddemo

Having specified all the relations among the \pp s of an
adjusted reducible homeomorphism $\phi_1$, the next step is to
to construct a homeomorphism isotopic to $\phi_1$ that eliminates as many
of these relations as possible. 
This is done by 
coalescing any \po s that are periodic Nielsen equivalent and
collapsing down those that are collapsible. In this process it may
be necessary to make minor changes in the topological type of
the components of $\phi_1$.

We first need a proposition about the existence of
homeomorphisms with specified fixed point behavior. The proof is
routine but we shall need some special features of the
homeomorphisms contained  in the proof. The statement
and use of Proposition 1.6 is similar to those of Theorem 4
in [Hp1]. The  {\de index} of a
fixed point $p$ with respect to $f$ is denoted $I(p, f)$ (see
for example, [J3]). The index of a  \pp\ $p$ with
$per(p, f) = n$ is $I(p, f^n)$.

\proclaim{Proposition 1.6:}  Let $M$ be a compact, connected
orientable surface with genus $g$ and $b$ boundary components
and let $n$ be a given positive integer. There exists a
homeomorphism $H:M\ra M$  isotopic to the identity which has $n$ fixed
points and no other \pp s. Further, except for the cases where
$\chi(M) = 0$ and $n=1$, the fixed points of $H$ have nonzero
index.
\endproclaim
	
\demo{Proof} The homeomorphism $H$ will be the
time-one map of a flow and is thus isotopic to the identity. The
first (and rather special) case is the sphere ($g=b=0$). When
$n=1$, let $H$ be the flow that has a single parabolic fixed point
at the North Pole. For $n=2$ use the North Pole/South Pole
gradient flow.  When $n=2k + 1$ or $n = 2k+2$, put k sink-saddle
pairs into the flows for $n=1$ and $n=2$, respectively.

The flow for the case $g=0$, $b=1$ and $n=2$ is shown in Figure 4(a).
 The dynamics of $H$ restricted to the boundary is rotation
by an irrational amount. The flow spirals outward from the
boundary towards the homoclinic loop.  The fixed point that is
not at the center of the bouquet of circles is a  spiral source
and the flow nearby spirals out towards the homoclinic loop.
The region  to the exterior of the bouquet consists of a single
parabolic component, \ie\ all the orbits in the region are
homoclinic to the fixed point at the center of the bouquet.  For
the case when $g=0$ and $b>1$, add more homoclinic loops with
boundary sources inside.  For other values of  $n$, add or
subtract homoclinic loops with spiral sources inside. The case
$g=0$, $b=1$ and $n=1$ is somewhat special in that the single
fixed point will be a saddlenode and have zero index.

\topinsert
\def\epsfsize#1#2{.7\hsize}
\centerline{\epsfbox{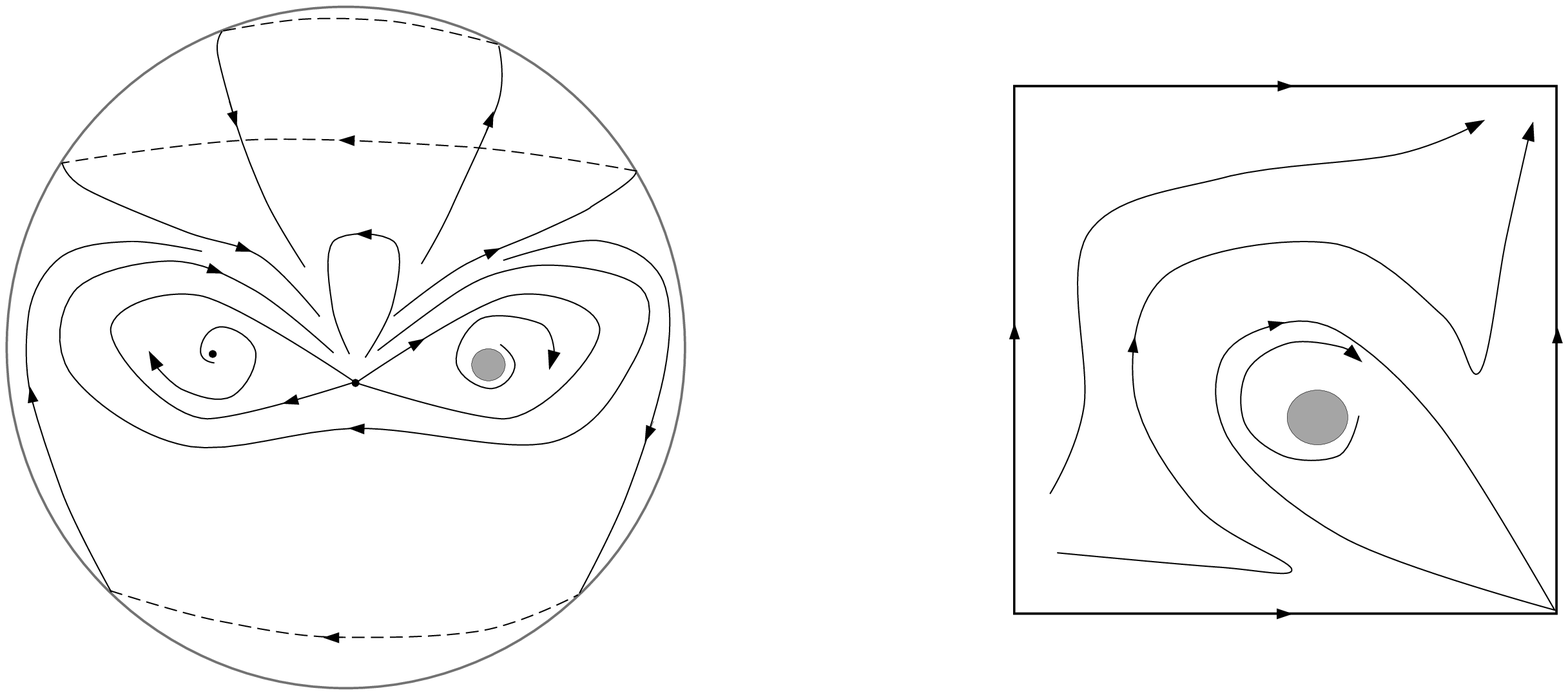}}
\botcaption{Figure 4} {{\bf (a)} A flow illustrating 
Proposition 1.6 for the case $g=0$, $b=1$, and $n=2$.
{\bf (b)} Same as (a) but for $g=b=n=1$.}
\endcaption
\endinsert

The flow for the case $g=b=n=1$  is shown in Figure 4(b).  The
arrows and labels on the edges indicate  the identification as
well as the flow direction. There is a fixed point at the
corners.  Once again, for larger values of $b$ and $n$ add more
homoclinic loops with boundary sources or spiral sources
inside.

For various values of $n$ when $g=1$ and $b=0$, 
use a similar construction but this time
there will be no boundary components inside homoclinic loops.
The case $n=1$ will have a fixed point with zero index.

A similar construction using a bouquet of circles also works
when
 $g>1$. An example is shown in Figure 5 for the case $g=4$,
$b=2$, and $n=3$. The letters on the boundary indicate identifications
\quad\qed
\enddemo

\topinsert
\def\epsfsize#1#2{.7\hsize}
\centerline{\epsfbox{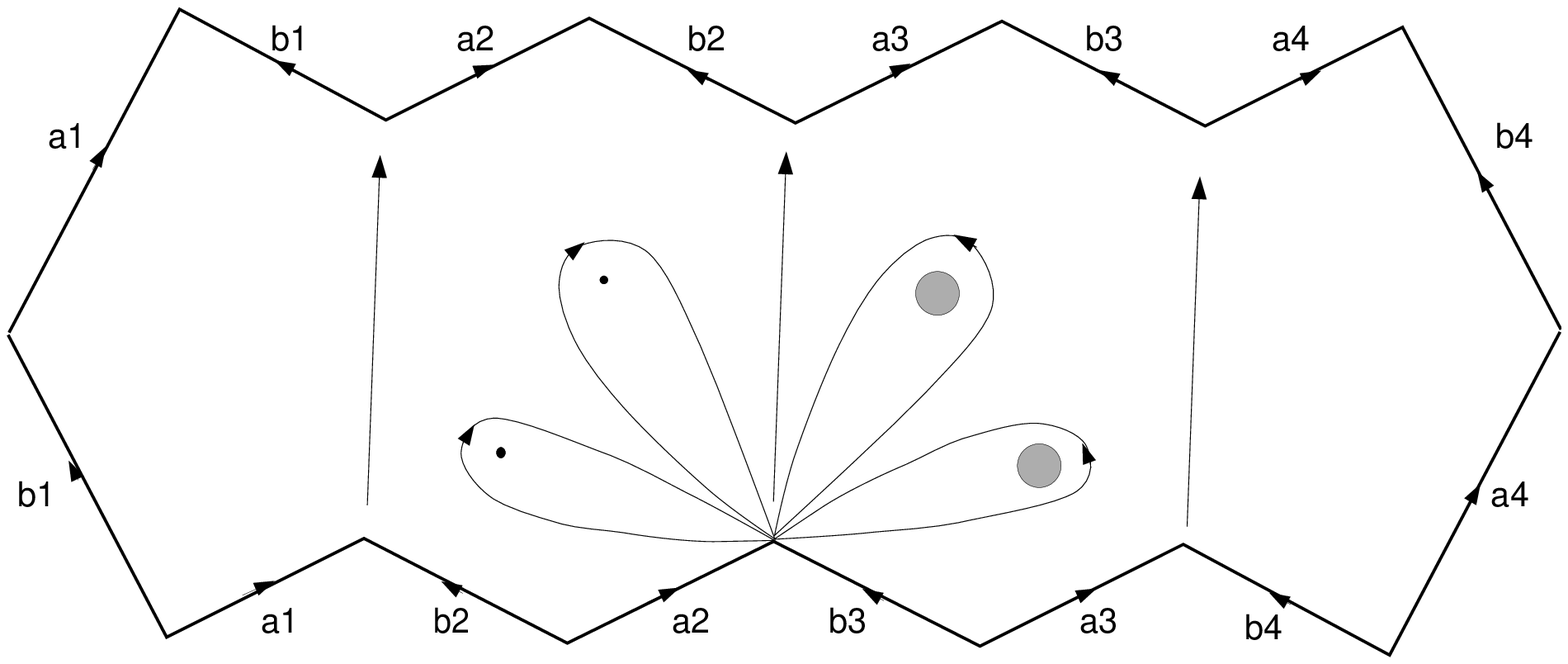}}
\botcaption{Figure 5}{A flow illustrating Proposition 1.6 
for the case $g=4$ and $b=n=2$.}
\endcaption
\endinsert

In a construction below the fixed points of the homeomorphism
$H$ from Proposition 1.6 will be lifted to periodic points in a
perhaps branched cover. For this application it is essential
that the lifted periodic points have nonzero index. If the
fixed point of $H$ is the image of regular points in the cover,
this is obvious. If the fixed point $p$ of $H$ is the image of  a
branch point, then the covering projection in a neighborhood of
$p$ looks like $z\mapsto z^k$ for some $k>1$.

When the dynamics in a neighborhood of a fixed point can be
decomposed into a finite number of sectors with each sector invariant
and $h$ of them
hyperbolic, $p$ of them parabolic, and the rest of them
attracting or repelling, then the index of the fixed point is $1+p/2-h/2$. 
If this fixed point is the image of a degree $k$ branch point, 
the index of the upstairs \pp\ under an appropriate iterate
that fixes it is $1+k(p/2 - h/2)$ if there is no local rotation
and $+1$ if there is. 
Thus, {\it even if the
index of the fixed point in the base is zero}, the index of the
periodic point in the cover is never zero. We note for future reference
that the  fixed points of the homeomorphisms constructed in Proposition 1.6 all
have a nice sector decomposition,  and so this fact is valid for these
maps.

We now begin the construction of the next refinement of the Thurston-Nielsen
canonical form. It alters an 
adjusted reducible homeomorphism $\phi_1$ to produce what will be
called a condensed homeomorphism. There are a number of cases to
consider corresponding to  different types of
$\phi_1$-components. In every case,   $N$ will  denote the
$\phi_1$-component under consideration  and $\phi$ denotes 
$\phi_1$ restricted to $N$. The end result in each case will
be a new homeomorphism $h$ to take the place of $\phi$  and a
new topological space $\Nc$ (perhaps not a manifold) to take the
place of $N$.

In a slight abuse of language, a boundary component of $N$ that
is the boundary circle of an untwisted reducing annulus will be
called an {\de untwisted boundary circle}. Let $U$ be the set of
all untwisted boundary circles in $N$.

The first four cases involve \fo\ components. In this case the
quotient space $N/\phi$ is denoted $N^\ast$. It is standard that
$\Na$ is a manifold and $\pi:N\ra \Na$ is a branched cover whose
upstairs branch points are precisely the branch \pp s of $\phi$.
In  the discussion of the various \fo\ cases Proposition 1.6
will be used to construct a homeomorphism $h^\ast:\Na\ra\Na$
with special properties.  The map $\ha$ is always isotopic to
the identity so it has a lift $h$ to $N$ that is isotopic to
$\phi$. Furthermore, $h$ will have the property that all its \pp
s have nonzero index.

{\bf Case 1:}  Assume that $\phi$ is \fo,  has no branch \pp s,
and $N$ has no untwisted boundary circles.  Pick a homeomorphism
$h^\ast:\Na\ra\Na$ that has one fixed point $p^\ast$ and no
other \pp s. Since in the case at hand $\pi : N\ra\Na$ is a
regular cover and $\chi(N)<0$, we have that $\chi(\Na) < 0$ and
so $I(p^\ast, h^\ast) < 0$.  Now lift $h^\ast$ to a
homeomorphism $h:N\ra N$ that has only one periodic orbit and
that orbit has nonzero index. In this case let $\Nc = N$.

{\bf Case 2:} Assume that $\phi$ is \fo, has no branch \pp s,
and $N$ has untwisted boundary circles. Let $U^\ast = \pi(U)$.
As in Case 1, $\pi : N\ra\Na$ is a regular cover, so $\chi(\Na)
< 0$. Using Proposition 1.6, choose $\ha$ which has a single
fixed point. Let $\Ra\subset\Na$ denote the subset of the
bouquet of circles used in the construction of $\ha$ which has
the property that each circle in $\Ra$ encloses an element
of $U^\ast$.  Now lift $h^\ast$ to a homeomorphism $h:N\ra N$
that has one periodic orbit and that periodic orbit has nonzero 
index. If $B^\ast$ is the open region
bounded by $\Ra$ and $U^\ast$, let $\Nc = N -
\pi^{-1}(B^\ast\cup U^\ast)$.

{\bf Case 3:} Assume that $\phi$ is \fo, has  branch \pp s, and
$N$ has no untwisted boundary circles. Let $P^\ast$ denote the
projection of the set of all branch \pp s to $\Na$.  Using
Proposition 1.6 construct  $\ha$  which has fixed points only at
points of $P\ast$ and no other periodic points. Let $h:N\ra N$
be the lift of $\ha$. By the remark after Proposition 1.6,  each
branch \pp\ in $N$ has nonzero  index under  $h$. In this case
let $\Nc = N$.

{\bf Case 4:} Assume that $\phi$ is \fo, has  branch \pp s, and
$N$ has untwisted boundary circles.  Using Proposition 1.6
  construct $\ha$ so that it has fixed points only at points of
$P^\ast$ and the points of $P$ have nonzero index under the lift
$h$. As in Case 2, let  $\Ra\subset\Na$ denote the subset of the
bouquet of circles used in construction of $\ha$ which has the
property that each circle in $\Ra$ encloses an element of
$U^\ast$. In
this case one may have to adjust $\Ra$ to make sure that no point 
of $P^\ast$ is
enclosed by a circle of $\Ra$.  By the remark after Proposition
1.6,   each branch \pp\ in $N$ has nonzero  index under  $h$.
Define $\Nc$ as in Case 2.

{\bf Case 5:} Assume $\phi$ is boundary-adjusted pA.  In this
case let $h = \phi$ and $\Nc = N$.

The next step in the process is to produce a map $\Phi$ that is
isotopic to the given adjusted reducible map $\phi_1$ by gluing
the new versions of the components and maps together.
There are once again a number of cases depending on the
character of the adjacent components. In every case the new map
$\Phi$ will be equal to the new component map $h$ on the
interior of the new component $\Nc$. In certain cases a closed
pA component or a reducing annulus will have its topological type altered by the
gluing. The structure of and the dynamics on  reducing annuli
must also be specified.

For adjacent $\phi_1$-components $N^1$ and $N^2$, we will denote
the corresponding new maps and spaces constructed in Cases 1--5
above as $h^1$, $N^{\circ 1}$, etc., and the old maps as
$\phi^1$ and $\phi^2$.  A reducing annuli between $N^1$ and
$N^2$ will be called $A$, and its $\phi_1$ orbit denoted by
$o(A)$.

If $N^1$ and $N^2$ are both pA and $A$ is twisted, no alteration
is necessary. This means that  $\Phi$ restricted to $o(A)$ is
the same as $\phi_1$ restricted to $o(A)$.

For pA components that adjoin an untwisted unflipped
annulus $A$, we may assume
that the dynamics on the two boundaries of $A$ are the same so we simply
glue these boundaries together eliminating the reducing
annulus.

If $N^1$ and $N^2$ are  pA and \fo, respectively, and  $A$ is
twisted then pick an isotopy of $\phi^2$ to $h^2$. Extend this
to an isotopy on $N^2\cup o(A)$ which starts at $\phi_1$
restricted to $N^2\cup o(A)$ and ends with a homeomorphism
called, say, $F$. Now compose $F$ restricted to $o(A)$ with a
push towards a boundary to insure it has no interior \pp s, and
let $\Phi$ restricted to $o(A)$ be the resulting map. The case
where both  $N^1$ and $N^2$ are \fo\ and $A$ is twisted is
similar.

Assume now that  $A$ is flipped and untwisted. 
As noted after Lemma 1.2, if $n$ is the 
least positive integer with $\phi^n(A) = A$, then 
$\phi^{2 n}$ has fixed points on both boundaries
of  $A$.  Further, $N^1$ and $N^2$ have the same type and by 
Lemma 1.3(d) they are pA. These imply that the boundary adjusted pA
on $N^1$ and $N^2$ has exactly one singularity
on each boundary of $A$. By the construction of the adjusted reducible
homeomorphism, $\phi^{2 n}$ restricted  to $A$ consists of the time $n$ map
of the flow shown in
Figure 3 composed with a rigid flip. Remove the region
between the boundary of $A$ and the homoclinic loops labeled $\Ra$
in Figure 3 and  glue the boundaries of $N^1$ and $N^2$ onto the
loops, identifying the singularities with the periodic point at
the ``pinch point''  of $\Ra$. A similar gluing is also
done on each component of the $\phi$ orbit of $N^1$, $N^2$, and
$A$. If prior to the gluing
the motion along the
homoclinic loops in $\Ra$ and the boundary of $N^1$ and $N^2$
has been properly adjusted,  then
one can define a homeomorphism on the $\phi$ orbit of the new space 
that is isotopic to $\phi$  and is equal to $\phi$ in the interiors
of the $\phi$ orbit of $N^1$,  $N^2$, and the new pinched version of $A$.

Since by Lemma 2(c) there are no untwisted reducing annuli
between adjacent \fo\ components, the last possibility is when
$N_1$ is pA and $N_2$ is \fo\ and $A$ is untwisted.  Now Cases
2 and  4 above were designed precisely so that in this
situation we could eliminate the reducing annuli and directly
glue the components together. If the \fo\ component is as in Case
4, this gluing will involve identifying certain peripheral \pp
s\ of the pA component  as they collapse down to
the branch \pp\ in the \fo\ component.
The only caution is this
process is that we must make sure that the direction of dynamics
on the bouquet of circles $R$ matches that on the boundary of
the pA component to which it is glued.  For this one may have to
adjust the choice of degenerate leaves that are collapsed when
constructing the boundary-adjusted pA map.

The homeomorphism  $\Phi$ that results from this process is
called {\de condensed}. Figure 6 illustrates an example of
the various stages of the construction. The center piece
shown in Figure 6(a) is a \fo\ component $N$ with  a 
component map which is simply a rotation by $180^\circ$ about a horizontal
axis. The fixed
point labeled $p$ is the only branch \pp. The adjacent components
are all pA, and all reducing annuli 
are untwisted. This means that the component $N$  falls into Case 4 above. 
The quotient
manifold $\Na$ and the bouquet of circles $\Ra$ is shown in Figure 6(b).
Figure 6(c) shows the manifold after the region in $N$ outside $R$
has been cut off and the components reglued.

\midinsert
\def\epsfsize#1#2{.7\hsize}
\centerline{\epsfbox{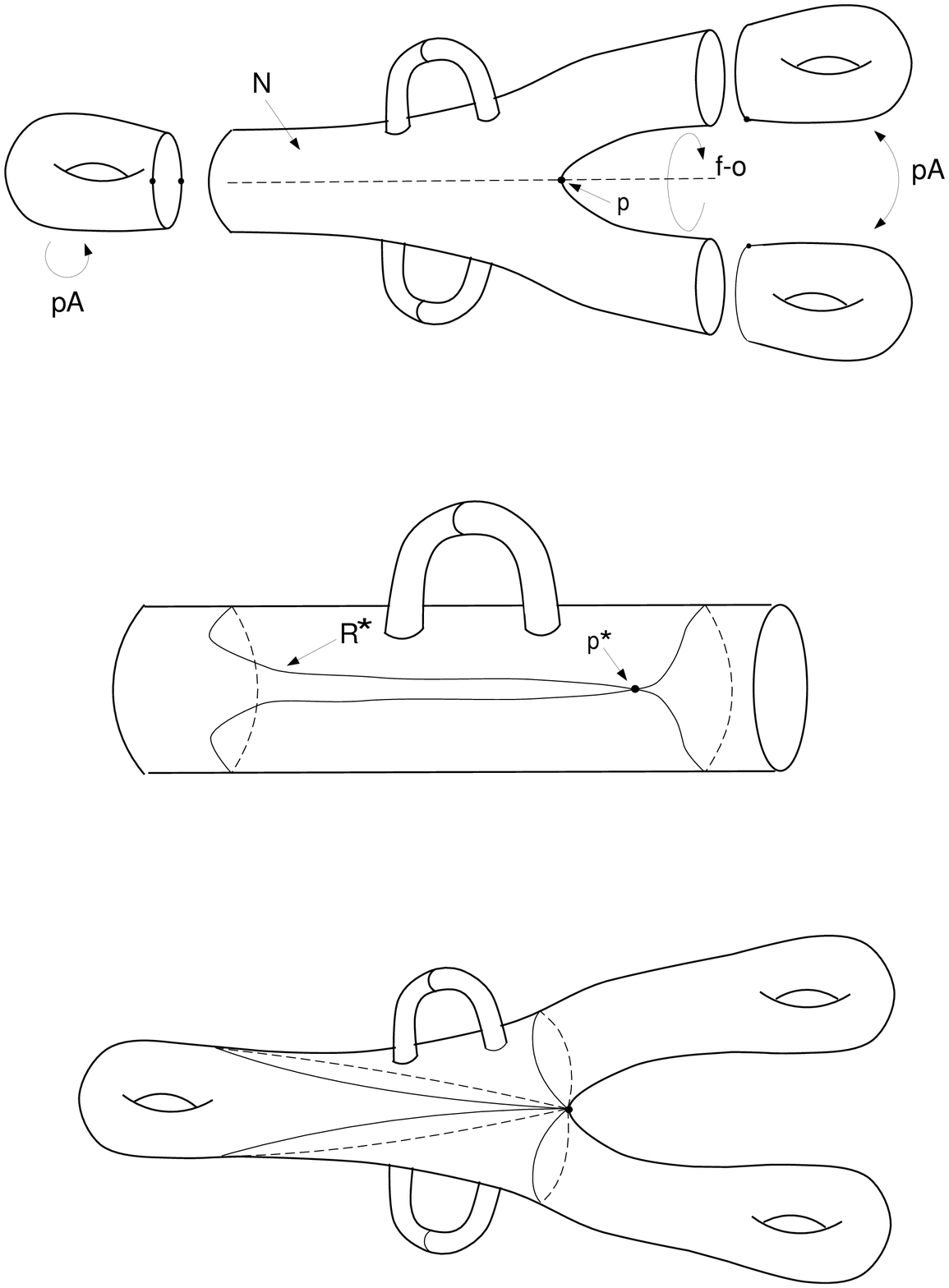}}
\botcaption{Figure 6}{{\bf (a)} A T-N reducible map on a 
genus $5$ surface.{\bf (b)} The quotient manifold for the central component
of (a) that is used in the Case 4 construction.
{\bf (c)} The reglued manifold with a condensed homeomorphism
in the isotopy class of (a).}
\endcaption
\endinsert

There was a certain amount of choice
involved in the construction of a condensed homeomorphism.
In particular, even up to conjugacy there is usually not a unique 
condensed homeomorphism in an isotopy class.
However, the construction was defined so that it eliminated most of the
relations among \pp s given in Proposition 1.5. The 
possible exceptions are peripheral and boundary pA \pp s and 
regular \pp s in \fo\ components that have no branch \pp s.
These \pp s all have the property
that they can be \PNE\ to another point on their orbit. It is immediate
from Lemma 2.1 that these \pp s are uncollapsible and therefore
unremovable. This implies  that these relations cannot be eliminated in
the isotopy class. 

Another important property of condensed homeomorphisms
is that each \pp\ has nonzero index. It is a standard fact
that this is true for the interior and boundary pA \pp s.
The construction guarantees this property for peripheral pA
and \fo\ \pp s as well as \pp s in the interior of flipped reducing
annuli.

Restricting attention to periodic {\it orbits\/} we have:

\proclaim{Theorem 1.7:} {\it  Every orientation-preserving
homeomorphism of an orientable  compact surface is isotopic to a
condensed homeomorphism. A condensed homeomorphism
has  no relations among its \po s and each \pp\ has nonzero index. }
\endproclaim

We have given the construction for this theorem in the case of
negative Euler characteristic. The other cases can be easily
handled case by case.

\head 2. Dynamically minimal models and persistence of \pp s\endhead
This section deals with the stability of \pp s under isotopy
in terms of \PNE\ and in terms of a stronger equivalence relation 
on \pp s called  \SNE.  
The main persistence result is Lemma 2.3 which builds on work of Jiang ([J3])
and Hall ([Hll]). This result allows us to conclude that all the
\pp s of a condensed homeomorphism persist under isotopy (Theorem 2.4).
Corollary 2.5 gives necessary and sufficient conditions for the
persistence of a \pp\ for a surface homeomorphism.

Although we focus here on homeomorphisms of surfaces, many of
the results in this section are true in a more general context.

Up to this point periodic Nielsen equivalence has been discussed
using arcs in the surface. There is an equivalent formulation using 
covering spaces.  
Let $\tM$ be the universal cover of $M$. Fix an
identification of $\pi_1(M) := \pi_1$ with the group of covering
translations of $\tM$. 
If we  fix a reference lift $\tf:\tM\ra\tM$ of $f$,
any lift of $f^n$ can be written as $\sigma \tf^n$ for some
$\sigma\in\pi_1$.
It is easy to check that two \pp s $x$ and $y$ are periodic Nielsen
equivalent if and only if there exists  lifts $\tx$, $\ty$ and an
element $\sigma\in\pi_1$  with $\sigma\tf^n(\tx) = \tx$ and
$\sigma\tf^n(\ty) = \ty$. Similarly, if $C$ is a simple closed curve
with $f^n(C) = C$, then $x$ and $C$ are periodic Nielsen
equivalent if and only if there exists  lifts $\tx$, $\tilde{C}$ and an
element $\sigma\in\pi_1$  with $\sigma\tf^n(\tx) = \tx$ and
$\sigma\tf^n(\tilde{C}) = \tilde{C}$ setwise. The equivalence class of 
$x$ under periodic Nielsen equivalence is $PNC(x, f)$.

Two lifts of $f$, $\tf_1$ and $\tf_2$, are said to be in the same
{\de lifting class} if there is an $\sigma\in \pi_1$ with $\sigma\tf_1
\sigma\inv = \tf _2$. The lifting class of $\tf$ is denoted $[\tf]$.
 If for a lift $\tx$ of $x$, the element
  $\sigma\in\pi_1$ is such that   $\sigma \tf^n(\tx) = \tx$, then
the {\de lifting class} of $x$ is $LC(x, f^n) =[\sigma \tf^n]$.
It is easy to check that this definition is independent of the
choice of the lift $\tx$ and if $(y,f) \pne (x,f)$, then $LC(y,
f^n) = LC(x, f^n)$.  It thus makes sense to talk about the lifting
class of a periodic Nielsen class.

A lifting class $\ell$ of $f^n$ is called {\de collapsible} if
there is a $\sigma\in\pi_1$ and an integer $k$ with $n = mk$ and
$1 \le k < n$ so that $\ell = [(\sigma\tf^k)^m]$. The next lemma
connects the notion of a collapsible lifting class with the notion
of collapsible \pp s introduced in the last section.  What is termed
``uncollapsible'' here is what was called ``irreducible''
in [J3]. page 65.  That terminology is not used to avoid confusion
with reducible maps.

\proclaim{Lemma 2.1} If $f:M\ra M$ is an orientation-preserving
 homeomorphism of a compact orientable $2$-manifold and $x\in
P_n(f)$, then following are equivalent:

{\leftskip=40pt\parindent=-18pt

(a) There exists a \pp\ $y$ so that $(x,f)\vdash (y,f)$.

(b) $LC(x,f^n)$ is collapsible.

(c) There exist integers  $k$ and $ m$ with  $1\le  k < n$ and
 $n = mk$, an element $\sigma\in\pi_1$, and a lift $\tx$ so that
$\tx$ is a \pp\ with period $m$  under $\sigma \tf^k$.

}
\endproclaim

\demo{Proof} First assume (a). If $per(x) = n$ and $per(y) = k$,
we have from the definition of collapsible  that $n = mk$ for some
$1<  m \le n$ and that for some $\sigma\in\pi_1$ $\sigma\tf^n(\tx)
= \tx$ and $\sigma\tf^n(\ty) = \ty$. But since $per(y) = k$, there
is also an $\alpha\in \pi_1$ with $\alpha\tf^k(\ty)= \ty$, and so
$(\alpha\tf^k)^m(\ty ) = \ty$.  But two lifts of $f^n$ with the
same fixed point must be the same and so $\sigma\tf^n =
(\alpha\tf^k)^m$. This implies (b).

Now assume (b). By assumption, if $\sigma\in\pi_1$ is such that
$\sigma\tf^n(\tx) =\tx$, then there are $\alpha, \beta \in \pi_1$
and $k$ with $1\le k < n$ so that $\sigma\tf^n = \alpha\inv
(\beta\tf^k)^m \alpha =(\alpha\inv \beta\tf^k\alpha)^m$. Thus $\tx$
is $m$- periodic under $\alpha\inv \beta\tf^k\alpha$ which is a
lift of $f^k$ and so (c) follows.

Finally assume (c). If $M$ has no boundary, the universal cover of
$M$ is either the plane or the sphere. The Brouwer Lemma states that
an orientation-preserving homeomorphism of the plane that has a
\pp\ also has a fixed point.
This implies that $\sigma \tf^k$ has a fixed point $\ty$. Thus
$(\sigma \tf^k)^m(\tx) = \tx$ and $(\sigma \tf^k)^m(\ty) = \ty$,
and so $x\vdash y$ which is (a). If $M$ has boundary then the
universal cover of $int(M)$ is the plane and thus if $x\in Int(M)$,
the proof is the same. If $x$ is on the boundary of $M$, collar
$M$ and extend $f$ so that it has no \pp s in the collar, and
proceed as in the case where $x\in Int(M)$.  
\quad\qed
\enddemo

It is important to note that the inclusion of condition
(a) in Lemma 2.1 is very dependent on the fact that we are in
dimension 2. The proof that (c) implies (a) uses the Brouwer Lemma 
which is not true in higher dimensions.

We now discuss a stronger  notion of equivalence of \pp s that was
first introduced in [AF].  It is convenient to include a parallel
discussion of periodic Nielsen equivalence. These two theories
share a general pattern in their development with any Nielsen-type
theory for periodic points.  We first introduce two notions associated
with isotopies.  A self-isotopy $f_t:f\simeq f$ is called {\de
contractible\/} if the corresponding closed loop in $Homeo(M)$ is
null-homotopic. An isotopy $f_t:f_0\simeq f_1$ is said to be a {\de
deformation\/} of a second isotopy $h_t:f_0\simeq f_1$ if the
corresponding arcs in $Homeo(M)$ are homotopic with fixed endpoints.

Assume that  $x_0\in P_n(f_0)$,  $x_1\in P_n(f_1)$ and $f_t:f_0\simeq
f_1$. The \pp s  $x_0$ and $x_1$ are {\de connected by the isotopy\/}
$f_t$ if there exists  an arc $\gamma:[0,1]\ra M$ with $\gamma(0)
= x_0$, $\gamma(1) = x_1$, and for all $t$, $\gamma(t)\in P_n(f_t)$.
If $(x_0, f_0)$ and $(x_1, f_1)$ are connected by some isotopy we
say that they are {\de connected by isotopy}.  Given a single map $f$
and  $x,y \in P_n(f)$, then $x$ is {\de strong Nielsen equivalent}
to $y$ (denoted $(x,f)\sne(y,f)$ or $x\sne y$) if $x$ and $y$ are
connected by a {\it contractible\/}
 isotopy $f_t:f\simeq f$.  The strong Nielsen class of a \pp\ $x$
is $SNC(x, f)$.  Two strong Nielsen classes are connected by an
isotopy if elements from each class are.

In the case of primary interest here ($M$ is a compact orientable
surface with negative Euler characteristic) all self-isotopies are
contractible ({\it cf.\/} [FLP], page 22).  However, the definition
of strong Nielsen equivalence is applicable in other situations so
the condition is explicitly mentioned  here.  The self-isotopy is
required to be contractible so that $x\sne y$ implies $x\pne y$.
This is a direct consequence of the fact that a contractible
self-isotopy always lifts to self-isotopy in the universal cover.
The fact that $x\sne y$ implies $x\pne y$ means that  a
periodic Nielsen class is composed of a disjoint collection of
strong Nielsen classes.  A strong Nielsen class is said to be {\de
collapsible} if its periodic Nielsen class is collapsible.

For fixed points there is no difference between the notions of
strong and periodic Nielsen equivalence ([J3], Theorem 2.13). 
The relationship of strong and periodic Nielsen equivalence 
for periodic orbits  is
clarified by considering the suspension flow of the given map.
Two periodic {\it orbits\/} are said to be  strong Nielsen equivalent
if  \pp s  from each orbit are.  It was shown in [J6] that two
periodic orbits are periodic Nielsen equivalent if and only if
their corresponding closed orbits are {\it homotopic\/} in the
suspension flow. One can show that for diffeomorphisms, two orbits
are strong Nielsen equivalent if and only if their corresponding
closed orbits are {\it isotopic\/} in the suspension flow.  This
indicates the close connection of strong Nielsen equivalence with
knot theory in dimension $3$ ({\it cf.\/} [BW]). It also indicates
that the distinction between strong and periodic Nielsen equivalence
vanishes in dimensions bigger than two. In this case the suspension
manifolds will be dimension $4$ or greater and in these dimensions
simple closed curves are homotopic if and only if they are isotopic.

One of the main purposes of a Nielsen-type theory for periodic
orbits is to find conditions that guarantee the persistence (in
the appropriate sense) of classes under isotopy (or homotopy). For
this the notion of index is crucial. Given a map $f$, any collection
of fixed points that is both open in $Fix(f)$ and closed in $M$
may be assigned an integer index ({\it cf.\/} [J3], page 11).

\proclaim{Lemma 2.2 \hbox{\rm (Essential Classes)}}
 An uncollapsible period $n$-periodic or -strong \hfill\break
Nielsen class is open in $Fix(f^n)$ and closed in $M$.
\endproclaim

\demo{Proof} We first prove the statement for periodic
Nielsen classes.  It is a standard fact that Nielsen fixed point
classes are closed in $M$ and open in the set of fixed points ([J3],
page 7). By applying this fact to the map $f^n$ the statement
follows after we have shown that when $PNC(x, f)$ is uncollapsible,
it is equal to the Nielsen class of $x$ under $f^n$. For this we
need only check that when $(x, f^n)\ne (y, f^n)$, then $y$ has
least period $n$. This is an immediate consequence of the equivalence
of Lemma 2.1 (a) and (b).

The statement for strong Nielsen classes follows from Lemma 4 in
[Hll] after we show that what is called
uncollapsible here implies what is called uncollapsible there.
Specifically, we must show that if $p\in P_n(f)$ is uncollapsible
in the sense used here, then the following condition holds: Whenever
there are sequences $g_j\ra g$ in $Homeo(M)$ and $q_j\ra q$ in $M$
with  $q_j\in P_n(g_j)$ and $(q_j, g_j)$ connected by isotopy to
$(p, f)$ for all $j$, then $q$ has least period $n$.

Fix lifts and $\sigma\in\pi_1$ so that $\sigma\tf^n(\tp) = \tp$.
Since each $(q_j, g_j)$ is connected by isotopy to $(p, f)$, for
all $j$ there exist lifts and equivariant isotopies $\tg_j\simeq\tf$
with $\sigma\tg^n_j(\tq_j) = \tq_j$ and $\tq_j\ra\tq$, which
implies $\sigma\tg^n(\tq) = \tq$ for some lift $\tg$ of $g$. Now
if $q$ had least period $k < n$, then there would be an element
$\alpha\in\pi_1$ with $\alpha\tg^k(\tq)=\tq$ and so $(\alpha\tg^k)^m
= \sigma\tg^n$ where $n = m k$.  This implies that $(\alpha\tf^k)^m
= \sigma\tf^n$, and so $p$ is collapsible, a contradiction. 
\quad\qed
\enddemo

A periodic or strong Nielsen class for which the index is defined
and is nonzero is called {\de essential}.

The next step in a Nielsen-type theory for periodic orbits is the
correspondence of \pp s in isotopic maps.  The notion  of connection
by isotopy is the type of correspondence appropriate to strong
Nielsen equivalence. It was defined above the last lemma.  If
$f_t:f_0\simeq f_1$ and $x_i\in P_n(f_i)$, then $PNC(x_0, f_0)$
{\de corresponds} to $PNC(x_1,f_1)$ under this isotopy if there
exists $\sigma\in\pi_1$ and lifts with $\sigma\tf_i^n(\tx_i) =
\tx_i$ where $\tf_t:\tf_0\simeq \tf_1$ with $\tf_t$ an equivariant
isotopy.  Equivalently, the \pp s correspond under the isotopy if
there is an arc $\gamma:[0,1]\ra M$ with $\gamma(0) = x_0$,
$\gamma(1) = x_1$, and the curve $f_t(\gamma(t))$ is homotopic to
$\gamma(t)$ with fixed endpoints.

\proclaim{Lemma 2.3 \hbox{\rm  (Correspondence and Persistence of Classes under
Isotopy)}} \hfill\break Assume that $f_t:f_0\simeq f_1$ and $x_i\in P_n(f_i)$,
for $i=1,2$.  \smallskip

(a) Periodic Nielsen Equivalence

{\leftskip=60pt\rightskip=16pt\parindent=-18pt

(i)  If  $PNC(x_0, f_0)$  is uncollapsible and corresponds to
$PNC(x_1,f_1)$ under the isotopy, then $PNC(x_1,f_1)$ is uncollapsible
and \hfill\break $I(PNC(x_0, f_0), f_0^n) = I(PNC(x_1, f_1), f_1^n)$.

(ii) If  $x_0$ is contained in a uncollapsible, essential periodic
Nielsen class, then there exists a $z\in P_n(f_1)$ that corresponds
to $x_0$ under the isotopy.

} \smallskip

(b) Strong Nielsen Equivalence

{\leftskip=60pt\rightskip=16pt\parindent=-18pt

(i)  If  $SNC(x_0, f_0)$  is uncollapsible and is connected  to
$SNC(x_1,f_1)$ by the isotopy, then  $SNC(x_1,f_1)$ is uncollapsible
and\hfill\break $I(SNC(x_0, f_0), f_0^n) = \allowbreak I(SNC(x_1, f_1), f_1^n)$.

(ii) If  $x_0$ is contained in a uncollapsible, essential strong
Nielsen class, then there exists a $z\in P_n(f_1)$ and an isotopy
$f_t^\prime:f_0\simeq f_1$ that is a deformation of $f_t$ so that
$z$ is connected  to $x_0$ by the isotopy $f_t^\prime$.

} 
\endproclaim

\demo{Proof} If $PNC(x_0, f_0)$ corresponds to $PNC(x_1, f_1)$,
then by definition there exists a $\sigma\in\pi_1$ and lifts so
that $\sigma\tf_i^n(\tx_i) = \tx_i$ and $\tf_0$ and $\tf_1$ are
equivariantly isotopic. Now if $PNC(x_1, f_1)$ is collapsible, then
there exists $\alpha, \beta\in \pi_1$ with $\sigma\tf_1^n = \alpha
(\beta\tf_1^k)^m\alpha^{-1}$ and so $\sigma\tf_0^n = \alpha
(\beta\tf_0^k)^m\alpha^{-1}$, so $x_0$ is collapsible.

The proof of Lemma 2.2 shows that when $x_0$ is uncollapsible its
periodic Nielsen class is equal to its Nielsen fixed point class
under $f^n$.  The equality of the indices and part (a)(ii)
then follow from [J3], Theorem I.4.5.

The assertions in (b) follow from [Hll] using the result concerning
uncollapsibility shown in the proof of Lemma 2.2.  
\quad\qed
\enddemo

It is convenient to have terminology to describe \pp s that behave
as in the (ii) parts in the above theorem.  A \pp\ $x\in P_n(f_0)$
is called {\de persistent\/} if for each given homeomorphism $f_1$
with $f_t:f_0\simeq f_1$, there is an $x_1\in P_n(f_1)$ which
corresponds to $x_0$ under the isotopy.  The \pp\ is called {\de
unremovable\/} if $x_0$ is connected to $x_1$ by an isotopy
$f_t^\prime:f_0\simeq f_1$, with the  isotopy $f_t^\prime$ a
deformation of $f_t$.

The next theorem  asserts  the existence of a minimal representative
with respect to periodic and strong Nielsen equivalence in the category
of orientation-preserving homeomorphisms of compact orientable
surfaces. The corresponding theorem for Nielsen fixed point classes
is outlined in [J2] and [I] and given in full detail in [JG].
In stating the theorem it will be convenient to work with equivalence
classes of both \pp s and \po s. The periodic (strong) Nielsen
class of a \po\ is simply the union of the classes for all the \pp
s in the orbit.

Theorem 1.7 says that any homeomorphism of the type considered here
is isotopic to a condensed homeomorphism. The condensed
homeomorphism has no relations
among its periodic {\it orbits\/} and each \pp\ has nonzero index.
The only relations among periodic {\it points\/} of the condensed
homeomorphism  were specified above that theorem. They consist of
uncollapsible \pp s that can be \PNE\ to other \pp s on the same
orbit but have no other relations. Therefore using Lemmas 2.1 and 2.3 we
have:  

\proclaim{Theorem 2.4 \hbox{\rm  (Dynamically Minimal Representative)}}
 Each orientation preserving homeomorphism  of a compact, orientable
$2$-manifold is isotopic to a condensed homeomorphism.
Each nonempty periodic (strong) Nielsen class of \pp s for a condensed
homeomorphism is uncollapsible and essential and is therefore 
persistent (unremovable).
Further, each  nonempty periodic  (strong) Nielsen class
of \po s contains exactly one element. 
\endproclaim

Informally, this theorem says that condensed
homeomorphisms are the least complicated dynamically in their
isotopy class. More formally, given an isotopy class $\Omega$,
for each $n$,  let $PON(\Omega, n)$ denote the number of uncollapsible,
essential period-$n$ periodic Nielsen classes of \po s for any map in
$\Omega$. Theorem 2.3 guarantees that this number is well-defined
and further, that it is a lower bound  for the number of
period-$n$ periodic orbits of any element in the isotopy class.
Theorem 2.4 asserts the existence of a map in the isotopy class
that achieves this lower bound for all $n$. 
A similar remark holds for strong Nielsen classes.

It is important to note there is no theorem of this type 
even for fixed point theories for certain homotopy classes of maps
on surfaces that are not homeomorphisms ([J4] and [J5]).

\proclaim{Corollary 2.5}  Let $f$ be an orientation-preserving
homeomorphism of a compact, orientable $2$-manifold. A \pp\ $x\in
P_n(f)$ is persistent (unremovable) if and only if its periodic
(strong) Nielsen class is uncollapsible and essential.
\endproclaim

\demo{Proof} Necessity follows from Lemma 2.3. Assume then
that $x$ is persistent. Let $\Phi$ be a condensed homeomorphism
that is isotopic to $f$. By definition of persistent, there is a
$z\in P_n(\Phi)$ that corresponds  to $x$ under the isotopy. By
Lemma 2.3(a), $PNC(x,f)$ is uncollapsible and essential because
$PNC(z,\phi)$ is. The proof for strong Nielsen equivalence is
virtually identical.  
\quad\qed
\enddemo

\remark{Remark 2.6}
The main goal in the construction of   condensed
homeomorphisms was to produce a map in the isotopy class having
the least number of \pp s. It will perhaps clarify the work of the
last section to describe  how and why a TN-reducible homeomorphism
fails to achieve this goal. 

As defined in Section 1, a TN-reducible map
can have  \pp s in the interior of reducing annuli. 
When these \pp s are contained in unflipped reducing annuli they
can always be removed by a strong push towards the boundary of the
annulus.  With the exception of a pair of \pp s, all the \pp s in the
interior of a flipped reducing annulus can be removed using Lemma 1.2(b).
Each uncollapsible strong Nielsen class that remains  is
essential.  This means that the only \pp s that are removable are
those that are collapsible. The collapsible periodic points in a
TN-reducible map are regular \pp s in \fo\ components that contain
branch \pp s,  peripheral pA \pp s in adjoining untwisted reducing
annuli, and peripheral pA \pp s on the boundary of a 
flipped, untwisted reducing annulus. 
These collapsible \po\ are removed using Cases 3 and 4 and
the gluing process from Section 1.

This leaves the following types of unremovable \pp s:  (1) Interior
pA \pp s. (2) Boundary pA \pp s.  (3) Peripheral pA \pp s that are
not adjacent via an untwisted reducing annulus to  \fo\ components
that have  branch \pp s  and are also not adjacent to an untwisted, 
flipped reducing annulus. (4) Finite order \pp s that are in a
component that has no branch \pp s. (5) Branch \pp\ in \fo\ components.
 (6) Periodic points in the interior of flipped reducing annuli.

Although none of these \pp s can be removed, the  map can be made
dynamically simpler by coalescing
 \pp s in the same  strong Nielsen class.  Points of type (1), 
(5), and (6)  are alone in their strong Nielsen classes so no coalescing is
needed for these points.  However, there are several
reductions that may be possible with points of type (2).  
In the standard model of
a pA map there are  pairs of \po s on the boundary that are strong
Nielsen equivalent.  These were coalesced in the corresponding
boundary-adjusted pA. Periodic points of type (2) can also be strong
Nielsen equivalent to other points from their orbit that lie on the same
boundary component. However, these \pp s are uncollapsible and
therefore coalescing to other points on the same orbit is
not possible. Similar remarks hold for \pp s of types (3) and (4).

Periodic orbits of types (2),  (3) and (4) 
can also be strong Nielsen equivalent to \pp s not on their orbit 
as described in Proposition
1.5 (a) and (c). These strong Nielsen equivalent \po s 
are coalesced in  the gluing process
in the construction of condensed homeomorphisms and the use of
Proposition 1.6, respectively. 
\endremark

\remark{Remark 2.7}
The conditions that guarantee unremovability of an $x\in
P_n(f)$ given in [AF] are that $SNC(x,f)$ is essential and that
the points\hfill\break $\{x, f(x), \dots, f^{n-1}(x) \}$ are in different
Nielsen classes as fixed points of $f^n$.  It is easy to check that
this last condition implies that  $SNC(x,f)$ is not collapsible.
These Asimov-Franks' conditions yield the unremovability of \pp s
of type (1), (5), and (6) in the previous remark. These types of \pp
s constitute ``most'' of the unremovable \pp s in ``most'' of the
isotopy classes.
 The other types of unremovable \pp s do not, in general, 
 satisfy these conditions.

The results of Hall in [Hll] strengthen the Asimov-Franks' result by
removing the need for differentiability of the homeomorphisms as
well as providing a lower level condition that ensures unremovability.
Hall's conditions also guarantees the unremovability of certain
finite collections of \pp s. Unremovability rel compact invariant
sets is considered in [BdHl].
\endremark

\remark{Remark 2.8}
There is an inherent awkwardness in dealing with
dynamically minimal models and equivalence classes of \pp s and
\po s.  This is reflected in the statements of Theorems 1.7 and 2.4.
As an example, let $\phi$ be a boundary-adjusted pA map on a surface 
$M$ that has
two boundary components that are mapped one to the other. Suppose further
that $\phi^2$ restricted to each of these circle has rotation number
$1/2$. This map certainly has the least number of \pp s among maps 
in its isotopy class. 
However, the boundary  \pp s are not alone in their periodic (or strong) 
Nielsen class as they are periodic (strong) Nielsen equivalent
to the other point on their orbit that
is on the same boundary circle. On the other hand, if one restricts
attention to equivalence of \po s, one has thrown out the more
detailed information about the structure of the orbit provided by
equivalence of \pp s.
\endremark

\remark{Remark 2.9}
As remarked above, periodic Nielsen classes are comprised
of the disjoint union of strong Nielsen classes.  The results of
this section show that  a persistent  periodic Nielsen class for
a surface homeomorphism contains exactly one unremovable strong
Nielsen class.
\endremark

\remark{Remark 2.10}
The conditions given in Corollary 2.5 that are necessary
and sufficient for unremovability have a heuristic interpretation
in terms of bifurcation theory. The intuitive idea is that a \pp\ is
removable exactly when there is a one-parameter family of homeomorphisms
with the property that the given \pp\ disappears in a bifurcation
at some point. This is essentially  the point of view taken in
[AF].  From this point of view the requirement that the strong
Nielsen class of the \pp\ be essential is necessary to guarantee
that the class of the \pp\ cannot disappear via a collection of
saddlenode bifurcations. The uncollapsibility condition prevents
disappearance via period-dividing bifurcations.
\endremark

\remark{Remark 2.11}
This section has presented a parallel development of two
Nielsen-type theories for \pp s. In such theories there are  other
components that have not been discussed here.  One such component
is the issue of coordinates for the equivalence classes. For periodic
Nielsen classes these coordinates are provided by lifting classes
 or by twisted conjugacy classes in $\pi_1$ ({\it cf.}
[J3] Section II.1). For maps isotopic to the identity, coordinates for
strong Nielsen classes are provided by the braid type which is the
conjugacy class in the mapping class group of the isotopy class of
the map on the complement of the orbit ({\it cf.} [Hll], Lemma 8 and
[Bd1]).

Another issue is how to compute the different classes and their
coordinates using algebraic or combinatoric information about the
maps. For surface homeomorphisms this is perhaps best accomplished
geometrically using train tracks ({\it cf.\/} [PH],   [BH],
[FM], [Ls], [Ky], and [Kl]).
Algebraic techniques are provided in [J3], [F], [FH], [HPY], and  [HY]

There are other Nielsen type theories for \pp s in addition to
those discussed here. The simplest way to describe these theories
is in terms of the suspension flow of the given homeomorphism. For
example, call two periodic orbits {\de Abelian Nielsen equivalent}
if the corresponding orbits are homologous in the suspension flow.
There is, of course, an equivalent formulation in terms of the
appropriate covering space and in terms of arcs in the surface.
For more information on such theories see, for example, [F]
or [HJ].
\endremark

\head 3. Persistence of pA orbits\endhead
The results of the  previous section show 
that a condensed homeomorphism is dynamically smallest in its
isotopy class and that its periodic orbits are present in
any isotopic map. If a condensed homeomorphism
has a pA component, the \po s are a small subset
of the interesting dynamics. This section concerns persistence
of these other orbits. 

Handel ([H]) and Fathi ([Ft]) both give persistence results for all the
dynamics of pA maps.  
Although global shadowing is not used explicitly, the point of view adopted 
here is closest to that of Handel (see Remark 3.3 below). 
The persistence of all the pA orbits from a condensed homeomorphism
is obtained by taking the closure of the persistent \po s.
This point of view is somewhat more natural here given the emphasis 
of the previous sections. It also avoids certain technical difficulties 
associated with the presence of a boundary on $M$ or
the presence of \fo\ components in a reducible map.

The first step is to use a given condensed homeomorphism $\Phi:M\ra M$
to define three pseudo-metrics on $\tM$, the universal cover of $M$ (in this 
paper {\it metric}  always refers to a topological metric, \ie, 
a map $M\times M\ra \reals$).
The pseudo-metrics are constructed using the invariant foliations
of the pA components of $\Phi$.  The measure attached to
an invariant foliation of a pA map assigns lengths to arcs
(\cf\ [FLP], expos\'e 5 \S II).  Given the lift of a pA component
$\tN_i$, the lift of the invariant unstable measured foliation of 
$\Phi$ restricted to $N_i$ gives a length $\ell_u^{(i)}(\gamma)$
to arcs $\gamma$ contained in $\tN_i$. For an arc $\gamma:[0,1]\ra\tM$
that intersects the lifts of pA components, $\tN_1, \dots, \tN_k$,
define $\ell_u(\gamma) = \sum_{i=1}^k \ell_u^{(i)}(\gamma\cap\tN_i)$.
Given $\tx, \ty\in\tM$, let 
$$\td_u(\tx,\ty)=\inf\{\ell_u(\gamma):\gamma  
\hbox{ is an arc connecting\ } \tx  \hbox{ and\ } \ty\}.$$
A pseudo-metric $\td_s$ is defined similarly using the lifts of 
stable foliations.
Let $\td_\Phi = \td_u + \td_s$. Note that these three pseudo-metrics
are equivariant, \ie\ for $\sigma\in\pi_1$, $d_s(\sigma\tx_0, \sigma\tx_1) = 
d_s(\tx_0, \tx_1)$, etc.

The pseudo-metric  $\td_\Phi$ projects to a pseudo-metric $d_\Phi$
on $M$ that is a metric when restricted to the interior of a pA
component ([FLP] pages 178--180). In the projected metric the distance 
between points in the same (connected) \fo\ component is zero.
Even in the case where there is only one component and
$\Phi$ is pA, $d_\Phi$ will not be an metric if $M$ has
boundary. This is  because  it assigns zero distance between 
pairs of points on the same boundary component (but these are the
only pairs of points that have zero separation). 
To  distinguish between these pseudo-metrics and an underlying
metric that makes $M$ a $2$-manifold (\eg\ a hyperbolic metric), the latter
will be called a {\de standard metric}. If $\rho$ is a metric or pseudo-metric
on $M$, the notation $(M,\rho)$ indicates the set $M$ with the topology
given by $\rho$.

Now  let $\lambda_\ast$ be the smallest
expansion constant among the pA components of $\Phi$. 
The crucial property of the pseudo-metrics is:
$$\eqalign{\td_u(\tPhi(\tx), \tPhi(\ty)) &\ge \lambda_\ast \td_u(\tx, \ty)\cr
\td_s(\tPhi^{-1}(\tx), \tPhi^{-1}(\ty)) 
&\ge \lambda_\ast \td_s(\tx, \ty).\cr}\eqno(\ast)$$
In particular, if $\tx$ and $\ty$ are a positive distance apart, then 
their separation as measured by $\td_\Phi$
grows exponentially under forward or backward iteration (or both).

We now return to the proofs promised in Section 1.

\demo{Proof of Lemma 1.1(b)} Recall that two  boundary
components $b_1$ and $b_2$ are related by a map $\phi$ 
if and only if there is an integer $n$, an element $\sigma\in\pi_1$, and 
lifts to the universal cover
$\tphi$ and $\tilde{b_i}$ so that $\sigma\tphi^n$ fixes the  $\tilde{b_i}$.
Now if $\phi$ is pA, property ($\ast$) and the fact that
$\td_\Phi$ is a metric on the interior of $\tM$ imply that $\sigma\tphi^n$
cannot fix two distinct boundary components of $\tM$. This shows that
distinct boundary components of $M$ cannot be related by $\phi$ 
and also yields Lemma 1.1(biv). The other parts of Lemma 1.1(b) follow
by considering points as well as boundary components.
\quad\qed
\enddemo

\demo{Proof of Lemma 1.4(a)}
By hypothesis, if $n = per(x, \phi_1)$, then there exists an
element $\sigma\in\pi_1$ and lifts with $\sigma\tphi^n(\tx) = \tx$,
$\sigma\tphi^n(\ty) = \ty$, and a $\tgamma$ that connect $\tx$ and
$\ty$.
Without loss of generality we may assume that $\Gamma_0$
is a simple closed geodesic with respect to a hyperbolic metric.
 We now identify $\tM$, the universal cover of $M$, with a subset
of the hyperbolic disk. Since $\gamma$ essentially intersects
$\Gamma_0$, there exists a lift $\tGamma_0$ so that
$\tgamma \cap \tGamma_0$ is nonempty and $\tx$ and $\ty$ are in
different components of $\tM - \tGamma_0$.

Since $\Gamma_0$ is a reducing curve, there is an $m > 0$ with
$\phi^m(\Gamma_0) = \Gamma_0$. This implies that
$(\sigma\tphi^n)^m(\tGamma_0)$ either equals $\tGamma_0$
or else is another (and therefore disjoint) lift of $\Gamma_0$.
But by compactness there are a finite number of lifts of
$\Gamma_0$ that intersect $\tgamma$ and so
$(\sigma\tphi^n)^m(\tGamma_0) = \tGamma_0$, which implies that
$\Gamma_0$ is related to both $x$ and $y$.

The proof of the second statement is similar.
\quad\qed
\enddemo

The next lemma contains the key idea of this section and
is taken almost directly from [H], Lemma 2.2. The 
difference here is in the use of the pseudo-metric $d_\Phi$ derived
from a reducible map  instead of the metric derived
from a pA map on a closed surface.
Extending a definition from Section 2,  if $f_0\simeq f_1$, 
the lifts $\tx_i$  of \pp s $(x_i,f_i)$ are said 
to correspond under the isotopy if there exists an equivariant 
isotopy $\tf_0\simeq\tf_1$ and an element $\sigma\in\pi_1$ with
$\sigma\tf_i^n(\tx_i)=\tx_i$, where $n = per(x_i, f_i)$.

\proclaim{Lemma 3.1} Let $f$ be an
orientation-preserving homeomorphism of a compact, orientable
$2$-manifold $M$ and $\Phi\simeq f$ be a condensed homeomorphism.
There exists a constant $C=C(f)$ such that whenever  $(x, \Phi)$
and $(y,f)$ are   \pp s with lifts $\tx$ and $\ty$ that correspond under
the isotopy, then $\td_\Phi(\tx, \ty) < C$.
\endproclaim

\demo{Proof} Fix equivariantly isotopic lifts $\tPhi$ and $\tf$
and  let 
$$R=\sup_{\tz\in\tM}\{\td_\Phi(\tPhi(\tz), \tf(\tz)), 
\td_\Phi(\tPhi^{-1}(\tz), \tf^{-1}(\tz))\}$$ 
and $C=2(R+1)/(\lambda_\ast-1)$.
Using the triangle equality and property ($\ast$) we have
$\td_u(\tPhi(\tx), \tf(\ty))\ge \lambda_\ast \td_u(\tx,\ty) - R$. Thus if
$\td_u(\tx,\ty)\ge C/2$,
then $\td_u(\tPhi(\tx), \tf(\ty))\ge 1 + \td_u(\tx,\ty)$, and
so  $\td_u(\tPhi^n(\tx), \tf^n(\ty))\ra\infty$ as $n\ra \infty$.

Similarly, if $\td_s(\tx,\ty)\ge C/2$, then
$\td_s(\tPhi^n(\tx), \tf^n(\ty))\ra\infty$ as $n\ra - \infty$.
Thus if $\td_\Phi(\tx,\ty)\ge C$, then $\td_\Phi(\tPhi^n(\tx), \tf^n(\ty))$
goes to infinity under forward or backward iteration (or both).

On the other hand, by hypothesis there exists an element $\sigma\in\pi_1$
with $\tPhi^N(\tx) = \sigma\inv\tx$ and $\tf^N(\ty) = \sigma\inv\ty$,
where $N=per(x,\Phi) = per(y, f)$. Since the pseudo-metric $\td_\Phi$ 
is equivariant this implies that
$\{\td_\Phi(\tPhi^n(\tx), \tf^n(\ty)): n\in{\Bbb Z}\}$ 
is finite, a contradiction.
\quad\qed
\enddemo

Note that a lemma of this type is certainly not true for
finite order maps using the lift of a standard metric.
The identity map gives a trivial example;
all fixed points are Nielsen equivalent, but there are lifts
that are fixed points arbitrarily far apart in the cover.
This is one reason the pseudo-metric $d_\Phi$ needs to vanish on \fo\
components.

The next theorem describes the persistence of the dynamics of a condensed
homeomorphism under isotopy. To get a statement that avoids 
pseudo-metrics it is necessary to pass to a quotient of $M$. 
Given a condensed homeomorphism
$\Phi:M\ra M$, define a equivalence relation $\bowtie$ on $M$ as follows.
If $z_0$ and $z_1$ are on the same boundary component of a pA
component, then $z_0 \bowtie z_1$. Boundary component here means
not only boundary components of the manifold itself, but also curves in
the interior of the manifold  which are the boundary to a pA component.
Points that are not on the boundary of a pA component are equivalent
only to themselves. Extend the relation $\bowtie$ so
that it is an equivalence relation. Note that 
the construction of a condensed homeomorphism can glue together
points from the boundaries of different pA components. These
points will be $\bowtie$-equivalent.

Let $M_p=M/\bowtie$. Fix a standard metric  on $M$ and 
give $M_p$ the metric induced by the projection $M\ra M_p$ (See Figure 7) 
Since $\Phi$ respects the relation $\bowtie$, it descends
to a homeomorphism $\Phi_p:M_p\ra M_p$.
In addition, since the projection $p$ is injective
on the interior of components of $\Phi$, we can continue to speak of
the components of the map $\Phi_p$.
If $\Phi$ is itself
pA, $M_p$ is just the closed manifold formed by collapsing
each the boundary components of $M$ to a point. In this case
the metric on  $M_p$ is equivalent to that induced by the
projection of $d_\Phi$ ([FLP], pages 178--180). In the general case,
on pA components in $M_p$ the given metric is equivalent
to the projection of $d_\Phi$. 

\topinsert
\def\epsfsize#1#2{.7\hsize}
\centerline{\epsfbox{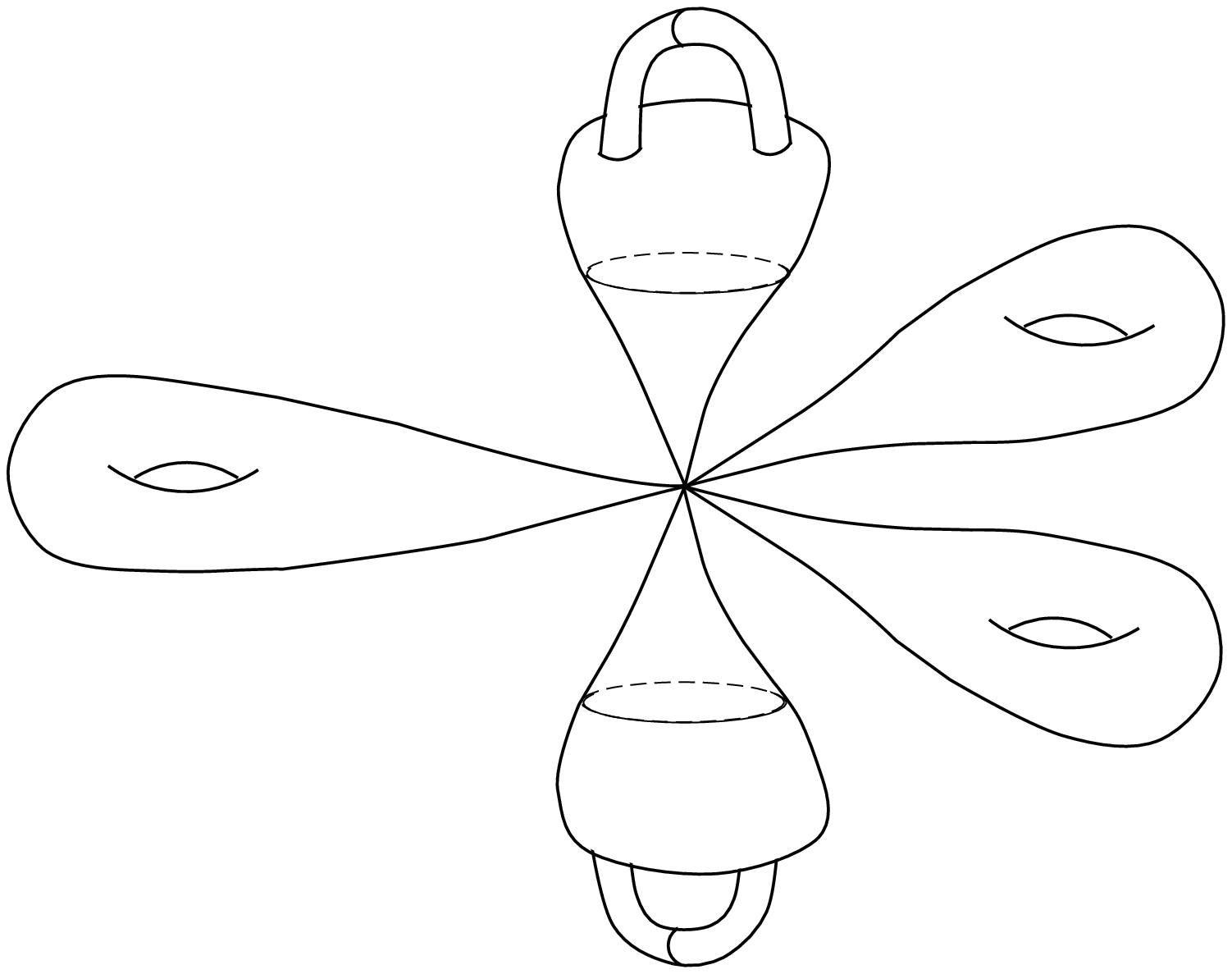}}
\botcaption{Figure 7}{The pinched manifold $M_p$ used in 
Theorem 3.2 derived from the condensed homeomorphism of Figure 6(c).}
\endcaption
\endinsert

\proclaim{Theorem 3.2} Let $f$ be an
orientation-preserving homeomorphism of a compact, orientable
$2$-manifold $M$ and $\Phi\simeq f$ be a condensed homeomorphism.
There exists a compact $f$-invariant set $\;Y\subset M$ and
a continuous map $\alpha:Y\ra M_p$ so that 
$\alpha \circ f_{| Y} = \Phi_p\circ \alpha$.
Further, $\alpha$ is homotopic to the inclusion and its
image contains all pA components and all \pp s of $\Phi_p$.
\endproclaim

\demo{Proof} 
Let $X_a$ and $X_f$ be the set of  \pp s of  $\Phi $
that are contained in the interior of a pA component or
the interior of a \fo\ component, respectively.
Pick a point $x$ from each \po\ of $\Phi$ contained in $ X_a\cup X_f$.
Let $\beta(x)$ be a \pp\ $(y,f)$ that is 
connected by isotopy to $(x,\Phi)$. By Theorem 2.4 and 
Lemma 2.3(bii)  such a point $y$
exists. Now extend $\beta$ by requiring that 
for $z\in X_a\cup X_f$, $\beta(z)$ is a \pp\ $(y,f)$ that is
connected by isotopy to $(z,\Phi)$, and further that 
$\beta\circ\Phi_{| X_a \cup X_f}=
f\circ\beta$.  By Theorem 1.7, $\beta$ is injective.  Let $Y_a=\beta(X_a)$,
 $Y_f=\beta(X_f)$,  and $\alpha_0 = \beta\inv$.

Since $d_\Phi$ is a metric on the interior of pA components
we may pick a standard metric $\rho$ on $M$ with $d_\Phi \le \rho$.
Let $\trho$ denote its lift to the universal cover. 
We first show that $\alpha_0:(Y_a\cup Y_f, \rho)\ra (M,d_\Phi)$ 
is uniformly continuous.
The proof depends on the following fact: If for $i = 1,2$, the \pp s
$(x_i, \Phi)$ and $(y_i, f)$ have lifts $\tx_i$ and $\ty_i$ that 
correspond under the isotopy, then for all $m$, 
$$\bigl| \td_\Phi(\tPhi^m(\tx_1), \tPhi^m(\tx_2)) - 
\td_\Phi(\tf^m(\ty_1), \tf^m(\ty_2))\bigr| \le 2 C.\eqno(\ast \ast)$$

To prove this fact note that 
$$\eqalign{\td_\Phi(\tf^m(\ty_1), \tf^m(\ty_2))
&\le \td_\Phi(\tf^m(\ty_1), \tPhi^m(\tx_1)) +
\td_\Phi(\tPhi^m(\tx_1), \tPhi^m(\tx_2))\cr
&\qquad+ \td_\Phi(\tPhi^m(\tx_2), \tf^m(\ty_2)).\cr}$$
By Lemma 3.1, the first and the
last terms on the right hand side are bounded by $C$. 
There is a similar inequality
obtained by switching all $\tPhi$'s and $\tf$'s and $\tx$'s and $\ty$'s.
These two inequalities together yield the fact.

Now assume to the contrary that $\alpha_0$ 
is not uniformly continuous. In this case
there exists some $\epsilon > 0$ such that for all positive $n$ there are
$y_1^{(n)}$ and $y_2^{(n)}$ with $\rho(y_1^{(n)}, y^{(n)}_2) < 1/n$ and
$d_\Phi(\alpha_0(y_1^{(n)}), \alpha_0(y_2^{(n)})) \ge \epsilon$.

Pick lifts $\ty^{(n)}_i$ of $y^{(n)}_i$ and $\tx_i^{(n)}$ of
$\alpha_0(y_i^{(n)})$  so that the inequalities of the previous paragraph
hold with the lifts in place of the points in the base. Fix an $M$ with
$\epsilon\lambda_\ast^M /2 > 4 C$. Since $\tf:(\tM,\trho)\ra (\tM,\trho)$
is uniformly continuous and $\td_\Phi \le \trho$  there is an $n$ so that
$$\max\{\td_\Phi(\tf^M(\ty^{(n)}_1), \tf^M(\ty^{(n)}_2)), 
\td_\Phi(\tf^{-M}(\ty^{(n)}_1), \tf^{-M}(\ty^{(n)}_2))\} < C.$$
By property ($\ast$), either
$\td_\Phi(\tPhi^M(\tx_1), \tPhi^M(\tx_2)) \ge\epsilon\lambda_\ast^M $
or 
$\td_\Phi(\tPhi^{-M}(\tx_1), \tPhi^{-M}(\tx_2)) \ge\epsilon\lambda_\ast^M $,
contradicting ($\ast \ast$) above.  Thus 
$\alpha_0:(Y_a\cup X_f, \rho)\ra (M,d_\Phi)$ is uniformly
continuous.

If $Y=Cl(Y_a)\cup Y_f$, where the closure is
taken with respect to topology given by $\rho$,
we can extend $\alpha_0$ to a continuous map, $\alpha_1:(Y,\rho)\ra(M,d_\Phi)$. 
The fact that  $M_p$ is obtained by identifying 
points on the boundaries of pA components whose $d_\Phi$ separation is zero
coupled with the fact that $X_f = \alpha(Y_f)$ is finite yields
that  $\alpha := p\circ\alpha_1:(Y,\rho)\ra M_p$ is continuous. 

Since \pp s of
$\Phi$ are dense in its pA components, the assertion about  the range of
$\alpha$ follows.
The fact that $\alpha \circ f_{| Y} = \Phi_p\circ \alpha$ is a straight
forward consequence of the definition of $\alpha$. If we pick a lift
$\tilde{Y}\subset\tM$, by virtue of Lemma 3.1 and  the fact that
$Y_f$ is finite,  $\alpha$ has a lift 
$\tilde{Y}\ra\tM$ that is a bounded distance from 
the inclusion. This implies that $\alpha$ is homotopic to the inclusion. 
\quad\qed
\enddemo

\remark{Remark 3.3}
 In [H] Handel proves Theorem 3.2 for pA maps on closed surfaces
using  global shadowing.
Given two isotopic maps $f$ and $\phi$ and equivariantly isotopic lifts
$\tf$ and $\tphi$, the pairs $(x,\phi)$ and $(y, f)$
{\de globally shadow} with respect to a metric
$d$ if there is a constant $K$ so that $\td(\tf^n(\ty), \tphi^n(\tx))<K$ for all
$n$, where as usual the tilde indicates lifts to the universal cover.

If $\phi$ is pA and  $Y_{gs}$ denotes the set of all points $(y, f)$ 
that globally shadow some point $(x,\phi)$, Handel shows that $f$ restricted
to $Y_{gs}$ is semiconjugate to $\phi$.  If the set constructed in the proof
of Theorem 3.2 is denoted $Y_{po}$, certainly $Y_{gs}$ can be larger 
than  $Y_{po}$.
A simple example is given by a DA map on the two torus (\cf\ [W], or the
appendix of [FR]). 
In this case $Y_{gs}$  is the entire torus  and $Y_{po}$ is just the basic set.
Clearly there are circumstances in which one or the other of
these sets would be most useful.

The metric used for global shadowing in a pA isotopy 
class on a closed manifold is 
$d_\phi$ as defined above.  As noted previously, 
when $M$ has boundary $d_\phi$ is only a pseudo-metric as it assigns zero 
distance to pairs of points on the same
boundary component. This means that with respect to $d_\phi$, a point 
$(y,f)$ globally shadows one point $(x, \phi)$ on the boundary if and
only if it globally shadows every point on the same boundary circle 
if and only if the orbit of its lift stays a bounded distance away 
from the orbit of the lift of the boundary circle. 

Using the pseudo-metric derived from a condensed homeomorphism  $\Phi$
one can prove a version of Theorem 3.2 via global shadowing.
In fact, most of the statements in [H] go through with minor changes. 
One again obtains a set $Y_{gs}$ that
has the pA components and \pp s of $\Phi$ as a factor.
There is, of course, still a difference between $Y_{gs}$ and $Y_{po}$.

As another example, let $\phi:M\ra M$ be pA and  assume that $M$ has boundary. 
Glue an annulus to a boundary component, extend $\phi$ 
in any manner, and call the new map $f$.
All the points in the new annulus under $f$ will globally shadow the 
boundary and therefore be in $Y_{gs}$, but only the points on the interior
boundary of the annulus will be in $Y_{po}$.
\endremark

\remark{Remark 3.4}
The intent in formulating Theorem 3.2 was to find
a model map that is a factor of every map in its isotopy class.
Unfortunately, to get such a result one is required to either 
use a degenerate topology coming from a pseudo-metric
(as in  $(Y,\rho)\ra (M, d_\Phi)$ is continuous) or else
use the pinched manifold $M_p$. The fact that some device
is necessary even in pA classes is illustrated by Figures 2(a) and 2(b).
These pictures
show possible boundary behavior for  two isotopic pA maps
that are conjugate on the interior of $M$. Any simple factor of these two maps
should clearly be the same on the interior of $M$, but there is no
boundary behavior that is a factor of the boundary dynamics of both maps.
The alternative adopted here was to collapse the boundary of the image
manifold to a point.
This yields a trivial point factor for the boundary dynamics.
\endremark

\refstyle{A}

\enddocument